\documentclass[12pt]{article} 
\usepackage[sectionbib]{natbib}
\usepackage{array,epsfig,fancyheadings,rotating}
\usepackage[]{hyperref}  
\usepackage{sectsty, secdot}
\usepackage{todonotes}
\sectionfont{\fontsize{12}{14pt plus.8pt minus .6pt}\selectfont}
\renewcommand{\theequation}{\thesection\arabic{equation}}
\subsectionfont{\fontsize{12}{14pt plus.8pt minus .6pt}\selectfont}

\usepackage[margin = 1in]{geometry}

\usepackage{amsmath}
\usepackage{amssymb}
\usepackage{amsfonts}
\usepackage{multirow}
\usepackage{amsthm}

\newtheorem{theorem}{Theorem}
\newtheorem{lemma}{Lemma}

\theoremstyle{definition}
\newtheorem{definition}{Definition}

\newtheorem{remark}{Remark}

\begin{document}

\markboth{\hfill{\footnotesize\rm INGRID DÆHLEN AND NILS LID HJORT} \hfill}
{\hfill {\footnotesize\rm MODEL ROBUST HYBRID LIKELIHOOD} \hfill}

\renewcommand{\thefootnote}{}
$\ $\par


\fontsize{12}{14pt plus.8pt minus .6pt}\selectfont \vspace{0.8pc}
\centerline{\large\bf Model robust hybrid likelihood }
\vspace{.4cm} 
\centerline{Ingrid Dæhlen and Nils Lid Hjort} 
\vspace{.4cm} 
\centerline{\it Department of Mathematics, University of Oslo}
 \vspace{.55cm} \fontsize{9}{11.5pt plus.8pt minus.6pt}\selectfont


\begin{quotation}
\noindent {\it Abstract:}
	The article concerns hybrid combinations of empirical and parametric likelihood functions. Combining the two allows classical parametric likelihood to be crucially modified via the nonparametric counterpart, making possible model misspecification less problematic. Limit theory for the maximum hybrid likelihood estimator is sorted out, also outside the parametric model conditions. Results include consistency of the estimated parameter in the parametric model towards a well-defined limit, as well as asymptotic normality after proper scaling and centring of the same quantity. Our results allow for the presence of plug-in parameters in the hybrid and empirical likelihood framework. Furthermore, the variance and mean squared error of these estimators are studied, with recipes for their estimation. The latter is used to define a focused information criterion, which can be used to choose how the parametric and empirical part of the hybrid combination should be balanced. This allows for hybrid models to be fitted in a context driven way, minimizing the estimated mean squared error for estimating any pre-specified  quantity of interest.

\vspace{9pt}
\par
\end{quotation}\par

\def\thefigure{\arabic{figure}}
\def\thetable{\arabic{table}}

\renewcommand{\theequation}{\thesection.\arabic{equation}}

\fontsize{12}{14pt plus.8pt minus .6pt}\selectfont

\section{Introduction and summary}
Parametric models are of course widely used in statistics. Procedures associated with these are often interpretable, and under model conditions they are in some sense asymptotically optimal. When the model fitted is incorrectly specified, however, these properties are no longer guaranteed to hold and instead bias can be introduced. Nonparametric methods, on the other hand, typically have little to no bias. They can, however, be quite variable, which might make the methods too imprecise to be useful in practice. In this article, we study hybrid likelihood theory, a middle ground between these two extreme approaches. We derive new and model agnostic results for the framework and greatly extend its applicability by proposing ways of deciding `how parametric' the method should be.

Hybrid likelihood theory was introduced in \citet*{hjort2018hybrid} and is a framework for combining the inference of nonparametric and parametric likelihood theory. For a given parametric model, with likelihood function $L_n(\theta)$, this is achieved by studying the hybrid likelihood function
\begin{align}\label{eq:HL def}
	\textup{HL}_n(\theta) = L_n(\theta)^{1-a}\textup{EL}_n[\mu(\theta)]^a,
\end{align}
where $\textup{EL}_n[\mu(\theta)]$ is the empirical likelihood function, descried in Section~\ref{chapter:EL}, a nonparametric counterpart to $L_n(\theta)$, evaluated at a context-driven control parameter $\mu(\theta)$, and $a$ is a tuning parameter in $[0,1)$ deciding how much weight should be put on the parametric and nonparametric parts of $\textup{HL}_n$. The control parameter $\mu$ is some pre-specified parameter we wish to estimate robustly, i.e.~without necessarily trusting the parametric model. Working with $\textup{HL}_n$ instead of either of its components separately allows us to choose a middle ground between purely parametric and nonparametric approaches. 

For the nonparametric part of $\textup{HL}_n$, the empirical likelihood function is used. This concept was first introduced in \citet{owen1988empirical}, and has since found numerous extensions and applications. For an overview, see e.g.  \citet{owen2001empirical} or \citet*{hjort2009extending}. So far, however, there have been few attempts to combine parametric and empirical likelihoods. Certain works study ways of combining the two methods when a parametric model is available for a subset of the data only, see e.g.~\citet{qin1994semi}. This is, however, quite different from the framework of  \citet{hjort2018hybrid}, where information from empirical and parametric likelihoods constructed using the whole sample are combined. Because of this, our hybrid likelihood framework turns out to be more related with semi-parametric methods like those of e.g.~\citet{olkin1987semiparametric} or \citet{mazo2021parametric}. Furthermore, the framework has certain conceptual similarities with methods developed from a rather different Bayesian semi-parametric setup in \citet{hjort1986discussion}. There a given parametric model has its likelihood modified by certain control sets, via a pinned down Dirichlet process on well-defined residuals. 

In this article we also derive new results for empirical likelihood theory. In Section~\ref{chapter:EL}, we use arguments similar to those of \citet*{molanes2009empirical} to derive an alternative characterization of the empirical likelihood function and the maximum empirical likelihood estimator (MEL estimator), defined and studied in \citet{qin1994empirical}. Rather than arguing as in that article, however, we use the framework of generalized estimating equations (GEE), see e.g.~\citet[chapter 5]{shao2003mathematical}. The method of GEEs allows us to study the limiting behaviour of the MEL estimator when plug-in estimators are used in the empirical likelihood function, i.e. $\hat\mu(\hat P)$, the maximizer of $\mu\mapsto\textup{EL}_n(\mu,\hat P)$ where $\hat P$ is some estimator. The behaviour of $\log\textup{EL}_n[\hat\mu(\hat P),\hat P]$ has been studied in \citet{hjort2009extending}, and \citet{molanes2009empirical} worked with $\hat\mu(\hat P)$ when the plug-in estimator is a function of $\mu$, equal to the maximizer of $\textup{EL}_n(\mu,P)$ for each fixed $\mu$. So far, however, no general result concerning the limiting behaviour of $\hat\mu(\hat P)$ for a general plug-in estimator has been derived. We believe our results fill this gap in the literature.

In Section~\ref{chapter:HL} we study the asymptotic properties of the hybrid likelihood function and its maximizer. \citet{hjort2018hybrid} did this, to some degree, by deriving asymptotic results for the hybrid likelihood function and its maximizer when the parametric model used in construction of the hybrid likelihood function is correctly specified, but their results do not hold true under model misspecification. Furthermore, we also derive results covering the case with plug-in estimators, which was not studied in \citet{hjort2018hybrid}. In Section~\ref{chapter:HL} we use the theorems of Section~\ref{chapter:EL} to show that the MHL estimator aims for the minimizer of a well-defined divergence $d_{a,\mu}$ and extend the results of \citet{hjort2018hybrid} to the case of possible model misspecification. We also mention that some preliminary results were reached in the master thesis \citet*{daehlen2021}. Our theorems also hold for the MEL estimator, making the theorems in Section \ref{section:MHLE} tools also for empirical likelihood theory in the context of possible model misspecification, i.e.~when $\textup{E}\,m(Y,\mu) = 0$ has no solution. This topic is apparently not previously studied in the literature. Both \citet{schennach2007point} and \citet{chib2018bayesian} study somewhat similar situations, but they both work with exponentially tilted empirical likelihoods, resulting in slightly different estimators than the MEL estimator.

For hybrid likelihood to be useful in practice, we need to choose how much weight should be put on the parametric and nonparametric part of $\textup{HL}_n$. In Section~\ref{section:The balance parameter}, we construct a focused information criterion (FIC) for models fitted with hybrid likelihood. The FIC was first introduced in \citet{claeskens2003focused} for models in a local misspecified modelling framework and in \citet{hjort2018hybrid} the authors derive a FIC similar in spirit to those of \citet{claeskens2003focused}. Our version of the FIC, however, is based upon the relatively newer version of the criterion introduced in \citet{jullum2017parametric}, with important modifications in \citet*{daehlenaccurate}. The criterion directly estimates the mean squared error of some pre-specified parameter of interest, which may or may not be equal to the control parameter. Using our proposed strategy to set $a$ should therefore, in theory, be asymptotically optimal and lead to more precise estimation than either of the methods can achieve on their own. This holds true without requiring model misspecification to disappear in the limit or to be absent altogether. Our derivations are based on the model robust results derived in Section~\ref{section:MHLE}, and would not be possible with previous methods like those of \citet{qin1994empirical}, \citet{white1982maximum}, or \citet{hjort2018hybrid}. We therefore believe we offer a new way to fit a hybrid parametric-nonparametric model in a way that is theoretically justified and guaranteed to lead to improved asymptotic inference.

After having discussed also regression settings, in Section~\ref{chapter:Regression} and presented some simulations in \ref{section:Simulations}, we go on to illustrate the hybrid likelihood machinery in Section~\ref{section:Wars}, pertaining to the important battle death series for all major interstate wars from 1823 to the 2022; see the Correlates of War data set, starting with \citet{CoW} and later curated by other scholars. Our analysis is a contribution to the ongoing discussion of whether wars have become somewhat less deadly over time.

\section{Empirical likelihood}\label{chapter:EL}
Before we can study the hybrid likelihood function and its maximizer, we define and take a closer look at the empirical likelihood function. This section will serve as motivation for the theory we derive for the hybrid likelihood function in Section \ref{chapter:HL}, in addition to establish notation which will be used throughout this article.

Empirical likelihood is a nonparametric tool for making inference about solutions  to equations on the form $0 = \textup{E}[m(Y,\mu)]$, where $Y\in\mathbb{R}^d$ is a random variable, $m\colon\mathbb{R}^{d+s}\to\mathbb{R}^q$ is some function and $\mu\in\mathbb{R}^s$ is the parameter we are solving the equation with respect to. Such expressions are called estimating equations, and $m$ is referred to as the estimating function. Many quantities can be expressed as solutions to equations on the form of $\textup{E}[m(Y,\mu)]=0$ for some function $m$. Examples include quantiles, moments, probabilities, and regression coefficients. 

Empirical likelihood theory is based on the study of the empirical likelihood function first introduced in \citet{owen1988empirical}. For an independent and identically distributed (i.i.d.) sample of random variables $Y_1,\ldots,Y_n\sim F$ in $\mathbb{R}^d$ and a fixed $m\colon\mathbb{R}^{d+s}\to\mathbb{R}^q$, the empirical likelihood function for the solution to $\textup{E}[m(Y,\mu)]=0$ takes the form
\begin{equation}\label{eq:Empirical likelihood function}
	\textup{EL}_n(\mu) = \max\left\{\prod_{i=1}^{n}w_i\,\colon\,\sum_{i=1}^{n}w_i = 1,\,\sum_{i=1}^nw_im(Y_i,\mu) = 0,\,w_i\geq0\right\}.
\end{equation}
The main result of empirical likelihood theory states that if $\mu_0$ is the solution to $\textup{E}[m(Y,\mu)]=0$ and the second moment of $m(Y,\mu_0)$ with $Y\sim F$ is finite, $-2\log\textup{EL}_n(\mu_0)$ converges in distribution to a chi-squared distributed variable with $q$ degrees of freedom. This result can be used to construct approximate non-parametric confidence intervals and hypothesis tests. For proofs, consult \citet{owen2001empirical} and \citet{hjort2018hybrid}. 

Using the Lagrange method, one can show
\begin{align}\label{eq:Empirical likelihood function lagrange}
	\textup{EL}_n(\mu) = \prod_{i=1}^n[1+\lambda_n(\mu)^Tm(Y_i,\mu)]^{-1},
\end{align}
where $\lambda_n(\mu)$ is the solution to
\begin{equation}\label{eq:Lagrange multiplier}
	0 = n^{-1}\sum_{i=1}^n\frac{m(Y_i,\mu)}{1+\lambda^Tm(Y_i,\mu)}
\end{equation}
under the condition that $1+\lambda^Tm(Y_i,\mu)\leq0$ for all $i=1,\ldots,n$. When no such $\lambda_n(\mu)$ exists, $\textup{EL}_n(\mu) = 0$. For derivations, see e.g.~\citet{owen2001empirical} or \cite{hjort2018hybrid}.

The right-hand side of \eqref{eq:Lagrange multiplier} is a mean of functions. Because of this, $\lambda_n(\mu)$ is a so-called Z-estimator for each fixed value of $\mu$. A Z-estimator is an estimator which can be expressed as the zero of an expression on the form $\sum_{i=1}^n\psi(Y_i,\theta)$ for some known function $\psi$. Such estimators have been well studied in the literature (see e.g.~\citet[Ch.~5]{van2000asymptotic}), and we will now use the relevant theory to derive an asymptotically equivalent expression for $\log\textup{EL}_n(\mu)$. For simplicity, we will leave $\mu$ out of the notation and replace $m(Y_i,\mu)$ by $M_i$ and $m(Y, \mu)$ for a general $Y\sim F$ by $M$.

Firstly, standard theory guarantees that a solution to \eqref{eq:Lagrange multiplier} exists with probability tending to one. In addition, the sequence of solutions has the following properties, 
\begin{equation}\label{eq:Limit of lagrange}
	\lambda_n\overset{\textup{pr}}{\to}\lambda\hspace{3mm}\text{and}\hspace{3mm}\sqrt{n}\left(\lambda_n-\lambda\right)=S^{-1}V_n+o_{\textup{pr}}(1)\overset{d}{\to}\text{N}(0,S^{-1}),
\end{equation}
where $\lambda$ is the solution to the equation $0 = \textup{E}[M/(1+\lambda^TM)]$ such that $0 = \Pr(1+\lambda^TM\leq0)$ and $S = \textup{E}\left[MM/(1+\lambda^TM)^2\right]$ and $V_n = n^{-1/2}\sum_{i=1}^nM_i/(1+\lambda^TM_i)$. For the above to hold, a relatively standard set of conditions needs to be satisfied, see e.g.~ \citet[Thms.~5.31 and 5.42]{van2000asymptotic}. For these results to be applicable, however, we need to assume there exists $L>0$ and an open ball, $\Lambda$, around $\lambda$ such that $\Pr(1+x^TM>L) = 1$ for all $x\in\Lambda$. This ensures that the function $(y,\lambda)\mapsto m(y,\lambda)/[1+\lambda^Tm(y,\lambda)]$ is smooth in $\lambda$ with probability 1 and that the theorems can be applied.

Since a solution to \eqref{eq:Lagrange multiplier} exists with probability tending to one, we can without loss of generality assume that $\textup{EL}_n$ can be written as in \eqref{eq:Empirical likelihood function lagrange}. Taking logarithms and Taylor expanding \eqref{eq:Empirical likelihood function}, reveals $\log\textup{EL}_n(\mu) + \sum_{i=1}^n\log(1+\lambda^TM_i) =-\sqrt{n}(\lambda_n-\lambda)^TV_n +\frac{1}{2}\sqrt{n}(\lambda_n-\lambda)^TS_n\sqrt{n}(\lambda_n-\lambda)+\epsilon_n$, where $S_n = n^{-1}\sum_{i=1}^nM_iM_i^T/(1+\lambda^T\mathit{M}_i)^2$. By \eqref{eq:Limit of lagrange} and the law of large numbers, this implies $\log\textup{EL}_n(\mu) + \sum_{i=1}^n\log(1+\lambda^TM_i)  = -\frac{1}{2}V_nS^{-1}V_n + \epsilon_n + o_{\textup{pr}}(1)$. Since open balls are convex, Taylor's theorem and the triangle inequality ensures $\left|\epsilon_n\right|\leq\lVert\lambda_n-\lambda\rVert^3\sum_{i=1}^n\lVert M_i\rVert^3/3L^3$ If the third moment of $\lVert M\rVert$ is finite, the law of large numbers and \eqref{eq:Limit of lagrange} guarantees that the right-hand side of this inequality is of order $O_{\textup{pr}}(1/\sqrt{n})$. In particular, this shows
\begin{equation}\label{eq:Pointwise limit of logEL}
	\log\textup{EL}_n(\mu) = -\sum_{i=1}^n\log[1+\lambda(\mu)^Tm(Y_i,\mu)] -\frac{1}{2} V_n(\mu)^TS(\mu)^{-1}V_n(\mu) + o_{\textup{pr}}(1),
\end{equation}
with probability tending to one. Here we have added $\mu$ to the notation for clarity and easier reference. 

Let $\mu_0$ solve the estimating equation $\textup{E}[M(\mu)] = 0$. Then $\textup{E}\{M(\mu_0)/[1+0\cdot M(\mu_0)]\}=0$. Hence, $\lambda(\mu_0) = 0$. Plugging this into \eqref{eq:Pointwise limit of logEL} shows $-2\log\textup{EL}_n(\mu_0) = V_n(\mu_0)^TS(\mu_0)^{-1}V_n(\mu_0)+o_{\textup{pr}}(1)\overset{d}{\to}\chi_q^2$ as $V_n(\mu_0)$ converges in distribution to a $N[0,S(\mu_0)]$ distribution by the central limit theorem. Hence, \eqref{eq:Pointwise limit of logEL} allows us to rediscover the main theorem of empirical likelihood theory. This is not the only application of the alternative characterization in \eqref{eq:Pointwise limit of logEL}, however. After dividing by $n$ we notice that the right-hand side of the equation is asymptotically equivalent to $-n^{-1}\sum_{i=1}^n\log[1+\lambda(\mu)^Tm(Y_i,\mu)]$. In the following section, we will use this fact to derive model robust theory in the context of hybrid likelihood. 

\section{Hybrid likelihood}\label{chapter:HL}
Hybrid likelihood inference was introduced in \citet{hjort2018hybrid} and is based on the study of the associated hybrid likelihood function. As explained in the introduction, this is a combination of the likelihood function of the data, $L_n$, and the empirical likelihood function of a pre-specified parameter of interest called the control parameter and denoted by $\mu$. We will assume that $\mu$ can be identified by the estimating equation $ \textup{E}[m(Y,\mu)] = 0$. Furthermore, for a given parametric model $f_\theta$ with $\theta\in\Theta\subseteq\mathbb{R}^p$, we will let $\mu(\theta)$ refer to the value the control parameter takes in the model $f_\theta$. With this notation, we  define the hybrid likelihood function as $\textup{HL}_n(\theta) = L_n(\theta)^{1-a}\textup{EL}_n[\mu(\theta)]^a$, where $\textup{EL}[\mu(\theta)]$ is the empirical likelihood function for $\mu$ evaluated at $\mu(\theta)$ and constructed using the estimating function $m$, and $a\in[0,1)$ is some tuning parameter deciding how much weight should be put on the parametric and nonparametric part of HL$_n$. The quantity $a$ is called the balance parameter.

The hybrid likelihood function can be interpreted as a more model robust version of the standard likelihood function, penalizing $\theta$s for which $\mu(\theta)$ is far from solving $0=n^{-1}\sum_{i=1}^nm(Y_i,\mu)$. If we, for instance, wish to fit a Weibull$(\lambda,k)$ model (parametrized such that the cumulative distribution function takes the form $F_{\lambda,k}(y) = 1-\exp(-(y/\lambda)^k)$) to some data while making sure that the mean is estimated robustly, we can use the hybrid likelihood function in place of the standard likelihood function. The mean satisfies $\textup{E}[m(Y,\mu)] = 0$ with $m(y,\mu) = y-\mu$. Hence, the empirical likelihood function can be constructed using this estimating function. Furthermore, in a Weibull$(\lambda,k)$ distribution, the mean takes the form $\lambda\Gamma(1+1/k)$ where $\Gamma$ is the Gamma function. Hence, in the hybrid likelihood function we would use $\mu(k,\lambda) = \lambda\Gamma(1+1/k)$. This procedure generalizes easily to higher moments, as we, in general, for the $\ell$-th moment can use $m(y,\mu) = y^\ell-\mu$ and $\mu(k,\lambda) = \lambda^\ell\Gamma(1+\ell/k)$. More generally, for a parametric model $f_{\theta}(y)$, the hybrid likelihood framework can ensure that estimation of $\textup{E}[h(Y)]$ is extra robust by using $m(y,\mu) = h(y) - \mu$ and letting $\mu(\theta)$ be the value of $\textup{E}[h(Y)]$ in the model $f_{\theta}$. A choice which might be of interest in many applications is to use $m(y,\mu) = I(y\in[b, c))-\mu$ for some pre-specified numbers $-\infty\leq b<c\leq\infty$ and $\mu(\theta) = F_\theta(c) - F_\theta(b)$ where $F_{\theta}$ is the cumulative distribution function in the model $f_{\theta}$. This would ensure that the fitted model estimates the probability of landing within $[b,c)$ extra robustly. This could be beneficial in situations where the fit in certain regions of the domain of $y$ are more important than others. 

In generalization of the above, we may construct a histogram controlled hybrid likelihood framework as follows. For bins $C_1,\ldots,C_k$, associated with a histogram, use estimating equations $m_j(u,\mu_j) = I(y\in C_j)-\mu_j$ for $j=1,\ldots,k-1$. This corresponds to using the probabilities $\mu_j=\Pr(Y\in C_j)$ for $j=1,\ldots,k$ as control parameters, and hence fit into the hybrid likelihood framework by using $m(y,\mu) =  (m_1(y,\mu),\ldots,m_{k-1}(y,\mu))^T$ and $\mu = (\mu_1,\ldots,\mu_{k-1})^T$ where $m_j(y,\mu) = I(y\in C_j) - \mu_j$ and $\mu_j(\theta) = \Pr_{\theta}(Y\in C_j)$ for $j=1,\ldots,k-1$. Here $\Pr_{\theta}(Y\in C_j)$ is the probability of $Y$ landing within $C_j$ in the parametric model $f_\theta$ for each $j$. By using these control parameters and corresponding estimating functions, the hybrid likelihood framework can be used to fit $f_\theta$ in a way ensuring that the probability of landing within each cell $C_1,\ldots,C_{k-1}$ is robustly estimated. This corresponds to ensuring that the histogram of the data with bins $C_1,\ldots,C_k$ fits well with the fitted model. 

\begin{figure}
	\centering
	\includegraphics[width=4in]{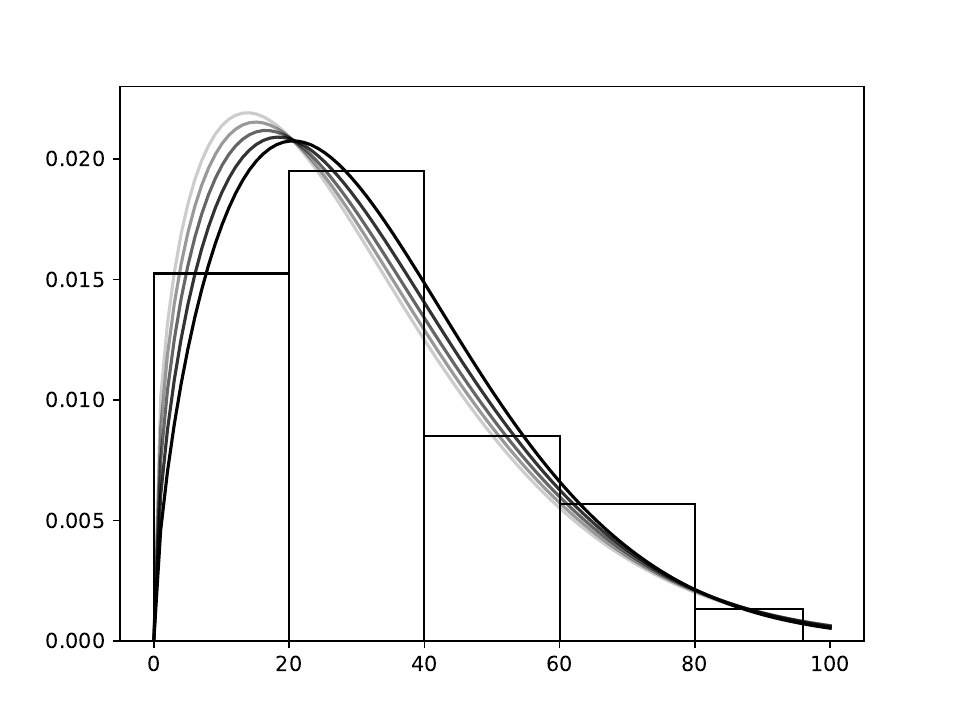}
	\caption{\label{fig:Weibull fit for each a} The figure shows a histogram of the life lengths for 141 men and women in Roman era Egypt together with the density of four fitted Weibull models. The model fitting was carried out by maximizing the hybrid log-likelihood function for $a=0,\,0.25,\,0.5,\,0.75,0.99$. The value of $a$ is indicated by the colour of the line with light values corresponding to low values of $a$.}
\end{figure}

We illustrate this idea on a data set taken from \citet{pearson1902v} containing life lengths of 141 men and women in Roman era Egypt. The data range from 1.5 to 96.0, and for this illustration we will fit a Weibull model to the life lengths. We will do this using the hybrid likelihood framework with $C_1 = [0,20), \ldots,C_5 = [80,100)$ and use $m$ as described above. In a In this case, $\mu_j(k,\lambda) = \exp[-(20(j-1)/\lambda)^k]-\exp[-(20j/\lambda)^k]$. The fit obtained using maximum likelihood is displayed in Figure \ref{fig:Weibull fit for each a}, together with the fits resulting from maximizing the hybrid likelihood function for different values of the balance parameter. From the figure we notice that as $a$ increases in value, the fitted density moves from the maximum likelihood estimate towards a curve which agrees somewhat more with the histogram of the data. This is a consequence of the empirical likelihood part of the hybrid likelihood function becoming more dominant as $a$ increases. As the Weibull-model fits the data quite well in this case, the fits are quite similar for different values of $a$. 

\begin{figure}
	\centering
	\includegraphics[width=4in]{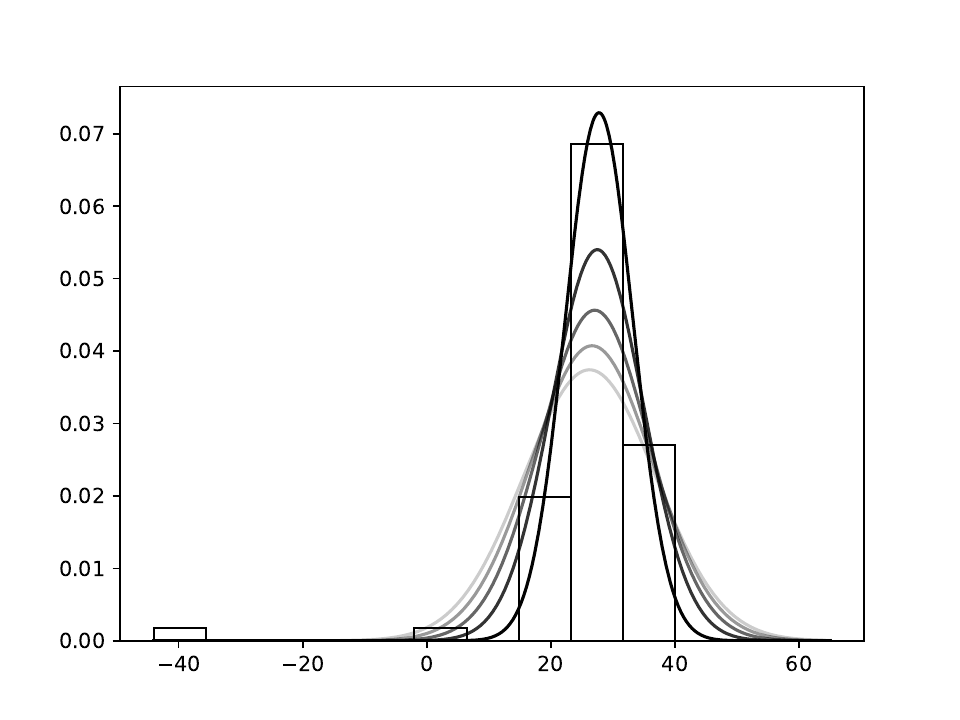}
	\caption{\label{fig:Normal fit for each a} The figure shows a histogram of transformed measurements of the speed of light together with the density of four fitted Normal models. The model fitting was carried out by maximizing the hybrid log-likelihood function for $a=0,\,0.25,\,0.5,\,0.75,0.99$. The value of $a$ is indicated by the colour of the line with light values corresponding to low values of $a$.}
\end{figure}

In other situations, however, the difference may be more substantial than for the data with Egyptian life lengths. For illustration of this, consider the dataset of 66 measurements of the speed of light, carried out by Simon Newcomb in 1882 (see e.g. Stigler, 1977). Each observation is measured in nanoseconds and recorded as the deviation from 24.8. The true value is 33.2.A histogram is shown in Figure~\ref{fig:Normal fit for each a} together with normal densities fitted using histogram controlled maximum hybrid likelihood for different values of $a$. From the figure we see that we get very different fits with different choices of balance parameter. This is a consequence of the outliers in the data set dominating the maximum likelihood fit and increased robustness as $a\to1$.

Now that we have introduced the framework of hybrid likelihood, we will study the asymptotic behaviour of the maximizer of the hybrid likelihood function. We will refer to this as the maximum hybrid likelihood estimator (the MHL estimator). \citet{hjort2018hybrid} consider the case when the parametric model $f_\theta$ is correctly specified. When this assumption holds true, the MHL estimator converges in probability to the true value of $\theta_0$ and at the speed of $O(1/\sqrt{n})$. As a result, the MHL estimator of $\mu$ also converges to the true value of $\mu$ in the underlying distribution of the data. Classic results concerning the empirical likelihood function in shrinking neighbourhoods of the true value (e.g.~\citet{hjort2009extending} and \citet{molanes2009empirical}) can therefore be used to derive limit results. When, on the other hand, the true underlying density is not assumed to be a member of the parametric family, there is no guarantee that the MHL estimator of the control parameter converges to the true solution of the estimating equation. Because of this, the results of e.g.~\citet{molanes2009empirical} and \citet{hjort2009extending} are not sufficient. Instead, we will use the arguments of the previous section.

\subsection{The hybrid likelihood function and its maximizer}\label{section:MHLE}
Similarly as we did for the empirical likelihood function, we start by deriving a pointwise asymptotically equivalent expression for $n^{-1}\log \textup{HL}_n(\theta)$, so $a\in[0,1)$ and let $h_n\colon\mathbb{R}^p\to\mathbb{R}$ be the hybrid log-likelihood function, defined as $h_n(\theta) = (1-a)\ell_n(\theta)+a\log\textup{EL}_n[\mu(\theta)]$, where $\ell_n$ is the log-likelihood of the data. Applying \eqref{eq:Pointwise limit of logEL} shows $h_n(\theta) = \sum_{i=1}\psi[Y_i,\mu(\theta)] + V_n[\mu(\theta)]^TS[\mu(\theta)]^{-1}V_n[\mu(\theta)] + o_{\textup{pr}}(1)$, where $\psi(y,\theta) = a\log f_{\theta}(y) - (1-a)\log\{1+\lambda[\mu(\theta)]^T$ $m[y,\mu(\theta)]\}$. In particular, the above shows that $n^{-1}h_n(\theta)$ is asymptotically equivalent to a mean. Because of this, we would expect its maximizer to be `almost' an M-estimator, a maximizer of a mean of functions (see e.g.~ \citet[Ch.~5]{van2000asymptotic} for definitions and properties of M-estimators). The following theorem shows that this intuition is indeed correct. The proof is given in the supplementary material, and to ease the notation, we will let $M_i(\theta)$ act as short-hand for $m[Y_i,\mu(\theta)]$ and $M(\theta)$ for $m[Y,\mu(\theta)]$ for a general $Y\sim F$. Furthermore, we will let  $\lambda(\theta)$ be the solution to $0  = \textup{E}\{M(\theta)/[1+\lambda^TM(\theta)]\}$, which is assumed to be unique. 

In the following theorem, we will need certain regularity conditions to be full-filled. Some of these are relatively straightforward, requiring the involved functions to be smooth. More specifically, for a certain set $\mathcal{N}$, we will need $\mu$ to be twice continuously differentiable and for $(\theta,\lambda)\mapsto \textup{E}\{\partial_\theta \psi(Y,\theta), M(\theta)^T/[1+\lambda^TM(\theta)]\}^T$ to be continuously differentiable. Furthermore, we will need both $m$ and $\log f_{\theta}$ to be regular on this set, in the sense that both the norm of $m$ and all partial derivatives of $m$ with respect to $\mu$ needs to be bounded by functions not depending on $\mu$ with finite fourth moments. Similarly, we will need all second order partial derivatives of $m$ as well as the score function in the model and all its partial derivatives to be bounded by square integrable functions. Under these conditions, the following theorem holds true.

\begin{theorem}\label{thm:HL limit marginal}
	Let $u$ be the score function in the model $f_\theta$ and fix $a\in [0,1]$. Suppose that $\hat\theta$ is the MHL estimator and unique critical point of $h_n$ in a open set $\Theta$. Define $\psi(y,\theta) = (1-a)\log f_\theta(y) - a\log[1+\lambda(\theta)^Tm(y,\theta)]$ and let $\theta_0$ be the maximizer and unique critical point of $\theta\mapsto \textup{E}\,\psi(Y,\theta)$ in $\Theta$. Then $\hat{\theta}$ is consistent for $\theta_0$ and
	\begin{align}\label{eq:Limit MHLE marginal}
		\hat{\theta} = \theta_0 + n^{-1}(J_a^*)^{-1}\sum_{i=1}^n\nabla_{\theta_0}\psi(Y_i,\theta_0) + o_{\textup{pr}}(1/\sqrt{n})
	\end{align}
	where $J_a^*= -H_{\theta_0}\textup{E}\,\psi(Y,\theta)$ and $\nabla_{\theta_0}$ and $H_{\theta_0}$ denotes the gradient and hessian matrix evaluated at $\theta_0$ respectively, provided there exists a neighbourhood $\mathcal{N}$ of $(\theta_0^T,\lambda(\theta_0)^T)^T$ on which the regularity conditions discussed in the paragraph preceding the theorem holds true, and in addition, that the function $(\theta,\lambda)\mapsto \textup{E}\{\partial_\theta \psi(Y,\theta), M(\theta)^T/[1+\lambda^TM(\theta)]\}^T$ has non-singular derivative at $(\theta_0^T,\lambda(\theta_0)^T)^T$, and that there exists some positive $L$ such that $\Pr(1+\lambda^TM(\theta)>L) = 1$ on $\mathcal{N}$. In particular, we have $\sqrt{n}(\hat{\theta}-\theta_0)\overset{d}{\to}\text{N}[0,(J_a^*)^{-1}K_a^*(J_a^*)^{-1}]$, where $K_a^* = \textup{Var}[\psi(Y,\theta_0)]$.
\end{theorem}
\begin{remark}\label{remark:Bounded support}
	Most of the assumptions of Theorem~\ref{thm:HL limit marginal} are more or less standard conditions ensuring sufficient regularity. The condition that, in $\mathcal{N}$, $\Pr(1 + \lambda^TM(\theta)>L) = 1$ for some $L>0$, however, warrants some discussion, as it excludes situations where $m(Y,\mu)$ has unbounded support, leading to several important situations not being covered. The perhaps most obvious example is $m(y,\mu) = y-\mu$ when $F$ has unbounded support. The assumption is, however, necessary, see e.g. \citet[Thm.~1]{schennach2007point}. In practice, this problem can be avoided by truncating data to lie within some arbitrary large and fixed interval.
\end{remark}

\begin{remark}
	For $a = 0$, $\theta_0$ is by definition equal to the minimizer of the Kullback-Leibler divergence from $f_{\theta}$ to $F$. Furthermore, we notice that $J_0^*$ and $K_0^*$ are equal to the fisher information matrix and the variance of the score function in the model $f_{\theta}$, respectively. Hence, Theorem \ref{thm:HL limit marginal} reduces to the results of \citet{white1982maximum} when $a = 0$. Additionally, even though the hybrid likelihood function is defined only for $a\in[0,1)$, the theorem above can be applied when $a=1$ and $h_n(\theta) = \log\textup{EL}_n[\mu(\theta)]$. In this case, it is most natural to let $\mu(\theta) = \theta$, for which Theorem \ref{thm:HL limit marginal} simplifies to a version of \citet[Thm.~1]{qin1994empirical}.
\end{remark}

If the true underlying distribution is a member of the parametric family fit to the data, there is a vector $\theta_0\in\mathbb{R}$ such that the distribution $F$ has density function $f_{\theta_0}$. In this case, $\mu(\theta_0)$ solves the estimating equation $\textup{E}\{m[Y,\mu(\theta_0)]\}=0$ and $\textup{E}\left\{m(Y,\mu)/\left[1+0\right]\right\} = 0$. Hence $\lambda[\mu(\theta_0)] = 0$, which together with direct computations shows that the derivative of $-\textup{E}\left\{\log\left[1+\lambda(\mu)^Tm(Y,\mu)\right]\right\}$ at $\mu(\theta_0)$ is zero. Because of this, $-\textup{E}\{\log\left[1+\lambda(\mu)^Tm(Y,\mu)\right]\}$ attains its maximum at $\mu(\theta_0)$, and as a result, both $\textup{E}[\log f_{\theta}(Y)]$ and $-\textup{E}\{\log[1 + \lambda(\theta)^TM(\theta)]\}$ are maximized by $\theta_0$, making this vector the maximizer of $\theta\mapsto(1-a)\textup{E}[\log f_{\theta}(Y)] - a\textup{E}\{\log[1 + \lambda(\theta)^TM(\theta)]\}$. Furthermore, direct computations show
\begin{equation*}
	\nabla_{\theta_0}\psi(y,\theta) = (1-a)u(y,\theta_0)  - a\textup{E}[ M'(\theta_0)]^T\textup{E}[M(\theta_0)M(\theta_0)^T]^{-1}\textup{E}[M'(\theta_0)],
\end{equation*}
Plugging this into Theorem~\ref{thm:HL limit marginal} shows $K_a^* = (1-a)^2J + a^2\xi_0^TW^{-1}\xi_0-a(1-a)(CW^{-1}\xi_0+\xi_0^TW^{-1}C^T)$, where $J$ is the Fisher matrix, $W$ is the covariance matrix of $m[Y,\mu(\theta_0)]$, $C$ is the covariance between $u(Y,\theta_0)$ and $m[Y,\mu(\theta_0)]$, and $\xi_0 = \textup{E}\,M'(\theta_0)$. In addition, $J_a^* = (1-a)J + a\xi_0^TW^{-1}\xi_0$ by the chain rule. Hence, if the model is correctly specified, Theorem~\ref{thm:HL limit marginal} simplifies to the corresponding results given in \citet{hjort2018hybrid}, making Theorem \ref{thm:HL limit marginal} a genuine model robust extension of the results of \citet{hjort2018hybrid}.

In certain situations, it can be beneficial to use a plug-in estimator inside the empirical likelihood part of $h_n$. Say that we, for instance, want to control for robust estimation of the variance $\sigma^2$ in a distribution. This can be achieved by using the hybrid likelihood framework with estimating function $m(y,\mu,\sigma^2) = (y-\mu,(y-\mu)^2-\sigma^2)$ this method however introduces the mean as a second control parameter. A perhaps more intuitive way could be to replace $\mu$ with the empirical mean of the observations $\hat{\mu}$ and use $m(y,\sigma^2) = (y-\hat{\mu})^2-\sigma^2$. The following result allows for such constructions. The proof is given in the supplementary material, and to improve readability, we adopt the following notation: Let $m\colon\mathbb{R}^{d+r+s}\to\mathbb{R}^q$ be an estimating function depending on a parameter $P$ which we wish to estimate using a plug-in estimator, i.e. $(y,\mu)\mapsto m(y,\mu,\hat P)$ is used as estimating function in the empirical likelihood function for some estimator $\hat P$. Let $\mu(\theta,P)$ be defined as the solution to $\textup{E}\{ m[Y, P,\mu(\theta,P)]\} = 0$ when $Y\sim f_{\theta}$. Lastly, let $\lambda(\theta,P)$ be the solution to $0 = \textup{E}\{ M(\theta,P)/[1+\lambda^TM(P,\theta)]\}$ for each $\theta$ and for each $P$. 

Again we will need a set of regularity conditions. These are relatively similar to those given in the paragraph preceeding Theorem~\ref{thm:HL limit marginal}. For a certain set $\mathcal{N}$, we will need $\mu$ to be twice continuously differentiable and for $(\mu,\lambda,P)\mapsto\textup{E}\{(\partial_\theta\psi_P(Y,\theta,P), M(\theta,P)^T/[1+\lambda^TM(\theta,P)], h(Y,P)^T)\}^T$ to be continuously differentiable. Furthermore, we will need the norm of $m$ and all partial derivatives of $m$ with respect to $\mu$ and $P$ to be bounded by functions not depending on $\mu$ or $P$ with finite fourth moments. Similarly, we will need all second order partial derivatives of $m$ as well as the score function in the model and all of its partial derivatives to be bounded by square integrable functions. Under these conditions, the following theorem holds true.

\begin{theorem}\label{thm:HL limit plug-in}
	Fix $a\in[0,1]$ and assume that $\hat{P}\in\mathbb{R}^r$ is an estimator satisfying $n^{-1}\sum_{i=1}^nh(Y_i,\hat{P}) = o_{p}(1\sqrt{n})$ for some function $h\colon\mathbb{R}^{d+r}\to\mathbb{R}^r$ and that $\hat P$ is used as a plug-in estimator for $P$ in the empirical likelihood part of $h_n$. Let $u$ be the score function in the parametric family, and $\hat{\theta}(\hat P)$ be the maximizer and unique critical point of $h_n$. Define $\psi_P(y,\theta,P) = (1-a)\log f_\theta(y) - a\log[1+\lambda(\theta,P)^Tm(y,\theta,P)]$. For each $P$ assume that $\theta(P)$ maximizes and is the unique critical point of $\theta\mapsto \textup{E}\,\psi_P(Y,\theta,P)$. Then the MHL estimator $\hat{\theta}(\hat P)$ is consistent for $\theta_0 = \theta(P_0)$ and
	\begin{align}\label{eq:Expansion for thetahat(phat)} 
		\hat{\theta}(\hat P) = \hat\theta(P_0) +\theta'(P_0) (\hat P-P_0) + o_{\textup{pr}}(1/\sqrt{n})
	\end{align}
	provided there exists a neighbourhood $\mathcal{N}$ of $\alpha_0 = (\theta_0^T,\lambda(\theta_0,P_0)^T,P_0^T)^T$ on which the regularity conditions stated in the paragraph preceeding the theorem hold true, and in addition that the function $(\mu,\lambda,P)\mapsto\textup{E}\{(\partial_\theta\psi_P(Y,\theta,P), M(\theta,P)^T/[1+\lambda^TM(\theta,P)], h(Y,P)^T)\}^T$  has non-singular derivative at $\alpha_0$. Furthermore, we need there to exists a positive $L$ such that $\Pr(1+\lambda^Tm(Y,\mu,P)>L) = 1$ on $\mathcal{N}$, and that on $\mathcal{N}$ the norms of $h$ and $\partial h/\partial P$ are bounded by $F$-square integrable functions not depending on $P$ or $\theta$.
\end{theorem}

\begin{remark}
	Although $\theta'(P_0)$ is well-defined, it can be hard to compute in practice as $\theta(P)$ is only implicitly defined. The implicit function theorem can, however, be used to find an explicit expression for the quantity.
\end{remark}

Let $V$ be the Jacobian matrix of $\theta\mapsto-\textup{E}\,h(Y,P)$ at $P_0$. Then we have $\hat P-P_0 = n^{-1}V^{-1}\sum_{i=1}^nh(Y_i,P_0) + o_p(1/\sqrt{n})$, by the conditions above and \citet[Thm.~5.21]{van2000asymptotic}. Furthermore, $\hat\theta(P_0)$ satisfies \eqref{eq:Limit MHLE marginal} with $m(y,\mu)$ replaced by $(y,\mu)\mapsto m(y,\mu,P_0)$ by Theorem \ref{thm:HL limit marginal}. Hence, the central limit theorem guarantees that $\sqrt{n}[\hat\theta(P_0) - \theta(P_0)]$ and $\sqrt{n}(\hat P -P_0)$ converge jointly in distribution to a central normal distribution with covariance matrix
\begin{equation*}
	\Sigma_\theta = \begin{pmatrix}
		(J_a^*)^{-1}K_a^*(J_a^*)^{-1},&(J_a^*)^{-1}C_\theta(V^{-1})^T\\
		V^{-1}C_\theta^T[(J_a^*)^{-1}]^T,&V^{-1}\text{Var}[h(Y,P_0)]V^{-1}
	\end{pmatrix}
\end{equation*}
where $V = \textup{E}[\partial_{P_0} h(Y,P)]$ and $C_\theta$ is the covariance between $\partial_{\theta_0}\psi_P(Y,\theta,P_0)$ and $h(Y,P_0)$ where $Y\sim F$. This together with \eqref{eq:Expansion for thetahat(phat)} ensures that $\sqrt{n}[\hat\theta(\hat P)-\theta(P_0)]$ converges to a normal distribution with mean zero and variance $(I_p,\mu'(P_0)^T)\Sigma_\theta(I_p,\mu'(P_0)^T)^T$.

\subsection{The hybrid divergence}\label{section:The limit function}
By Theorem~\ref{thm:HL limit marginal}, the maximum hybrid likelihood estimator converges in probability to the maximizer of $\textup{E}\,\psi(Y,\theta) = (1-a)\textup{E}[\log f_{\theta}(Y)] - a\textup{E}\{\log[1 + \lambda(\theta)^TM(\theta)]\}$. This function is useful for showing consistency of the MHL estimator, but the expression is somewhat hard to interpret, and it is not clear from the formula why we should want to maximize it. In this section, we will take a closer look at $\textup{E}\,\psi(Y,\theta)$ and show that maximizing this expression is equivalent to minimizing a divergence taking both model fit and precise estimation of the control parameter into account.

We start by considering $-\textup{E}\{\log[1 + \lambda(\theta)^TM(\theta)]\}$. Using the Lagrange method of multipliers, one can show that for each fixed $\mu$, the expression is the solution to the following optimization problem:
\begin{align}\label{eq:Divergence limit EL 1}
	\underset{w\colon\mathbb{R}^d\to\mathbb{R}}{\max}\,\left\{\textup{E}[\log w(Y)]\,\colon\,\textup{E}[w(Y)M(\mu)] = 0\,\text{and}\,\textup{E}[w(Y)] = 1\right\}.
\end{align}
Assume now that the true distribution of $Y$ has a density (pdf) or point mass function (pmf) $g$ and let $h(y) = w(y)g(y)$. Then solving \eqref{eq:Divergence limit EL 1} is equivalent to minimizing $\text{KL}(g,h)$ with respect to $h$ under the condition that $h$ is a pmf/pdf and $\textup{E}_h[M(\mu)]=0$. Here KL$(g,h)$ is the Kullback-Leibler divergence from $h$ to $g$ defined as $\textup{E}\{\log [h(Y)/g(Y)]\}$. Since subtracting a constant does not change the extrema of an expression, the above shows that the maximizer of $\textup{E}\,\psi(Y,\theta)$, i.e.~the limit of the MHL estimator, is also the minimizer of $\theta\mapsto d_a(g,f_\theta)$ where 
\begin{align}\label{eq:Distance HL}
	d_a(g,h) = (1-a)\text{KL}(g,h) + a\,\underset{h}{\min}\,\left\{\text{KL}(g,h)\,\colon\,\textup{E}_{h}\{M[\mu(h)]\} = 0\right\},
\end{align}
with $\mu(h)$ being the solution to the estimating equation $\textup{E}_h[M(\mu)] = 0$. The value of KL$(g,h)>0$ is zero when $g\neq h$. Since also the second term is non-negative by the previous arguments, $d_a(g,h)>0$ for all $h\neq g$. Furthermore, KL$(g,g) = 0$ and $\textup{E}_g\{M[\mu(g)]\} = 0$. Hence, $d_a(g,h) = 0$ only when $h = g$. Because of this, $d_a$ is a proper divergence, measuring the `distance' from the true density $g$ to $h$.

The first term in \eqref{eq:Distance HL} is small when the overall model fit is good, i.e.~when $h$ is close to $g$ overall. The second term, on the other hand, does not emphasise overall model fit, but is small when $\mu(h)$ is close to $\mu(g)$. Hence, $d_a$ measures the distance from $g$ to $h$ by a convex combination of the lack of overall model fit and error involved when using $\mu(h)$ in place of $\mu(g)$. By Theorem~\ref{thm:HL limit marginal}, the MHL estimator converges to $\theta_0$ such that $f_{\theta_0}$ minimizes $d_a$ in the class $f_\theta$ for $\theta\in\Theta$. Hence, the limit of the MHL estimator is the value of $\theta$ making $f_{\theta}$ as close as possible to the true density $g$, and distance is measured as a trade-off between lack of model fit and bias of the control parameter, and how much weight should be put on each term in the divergence is determined by the value of the tuning parameter. 

\section{The balance parameter}\label{section:The balance parameter}
To use the hybrid likelihood theory in practice, we need to choose a control parameter, $\mu$, and set the value of the balance parameter, $a$. The choice of the former should reflect what parameters are of particular interest in the specific application, and should not be chosen in an automatic way. In some situations, the same might be true for $a$. Usually, however, how much weight should be put on the parametric and nonparametric part of $\textup{HL}_n$ is dependent on the quality of the model fit and estimator of the control parameter. In this section, we will propose a way to choose $a$ taking these considerations into account.

\subsection{FIC for models fitted with hybrid likelihood}
In certain applications, estimation of a single parameter is the main point of the analysis. In such cases, models resulting in precise estimators of this parameter of interest are more desirable than those achieving a good overall fit. This is the main idea behind the focused information criterion. 

The focused information criterion (see \citet{claeskens2008model} for an introduction) ranks models by the quality of their estimate of the focus parameter. The information criterion is created to choose between different parametric models, but the estimators defined can in principle be used for other purposes. We propose to use the information criterion to choose the optimal value for the balance parameter $a$ used in construction of the hybrid likelihood function. This idea was briefly visited in \citet{hjort2018hybrid}, where the authors derive formulas for a focused information criterion for models in shrinking neighbourhoods of the true underlying density. More precisely, \citet{hjort2018hybrid} formulate an FIC for models on the form $f(y,\theta)$ when the true underlying distribution of the data is $f_{\textup{true}}(y,\theta,\gamma_0 + \delta/\sqrt{n})$ and $f_{\textup{true}}(y,\theta,\gamma_0) = f(y,\theta)$. Their version of the information criterion therefore requires the model to be asymptotically correctly specified, and to be able to compute their FIC one needs to define a wider model $f_{\textup{true}}(y,\theta,\gamma)$, which is then assumed to be correct. This is not too problematic in e.g.~the context of variable selection, but can be quite restrictive  for more general use. We will therefore instead create a fully model agnostic FIC, based upon the formulas and results of \citet{jullum2017parametric} with important modifications in \citet{daehlenaccurate}. These articles find asymptotically unbiased approximations to the mean squared error of estimators of the focus parameter in a large number of cases. Their results can be used to find estimators of the mean square error of maximum hybrid likelihood estimators of focus parameters. In the following theorem, we state a version of the results of \citet{jullum2017parametric} and \citet{daehlenaccurate}  valid in the hybrid likelihood context. 

\begin{theorem}\label{thm:FIC theorem}
	Fix a balance parameter $a$ and let $\hat\theta$ be the MHL estimator. Let $\xi\in\mathbb{R}$ be some focus parameter and $\hat{\xi}$ a consistent estimator of $\xi$ satisfying
	\begin{align}\label{eq:Influence function FIC}
		\hat{\xi} = \xi + n^{-1}\sum_{i=1}^n\phi(Y_i,\xi) + \epsilon_n
	\end{align}
	for some function $\phi\colon\mathbb{R}^d\to\mathbb{R}$ and $\epsilon_n = O_{\textup{pr}}(1/n)$. Assume further that all conditions of Theorem~\ref{thm:HL limit marginal} hold true and that there is a continuously differentiable function $g\colon\mathbb{R}^p\to\mathbb{R}$ such that $\xi = g(\theta)$ in the model $f_\theta$. Then $\hat{b}^2 - n^{-1}(2bc + \kappa - \tau)$  is only $o(1/n)$ away from the true MSE of $g(\hat{\theta}) $ in expected value where $\hat{b} = g(\hat{\theta})-\hat{\xi}$, $c = \lim_{n\to\infty}n\textup{E}(\epsilon_n)$, $\kappa = \textup{Var}[\nabla g(\theta_0)^T(J_a^*)^{-1}\psi_{\theta}(Y,\theta_0) - \phi(Y,\xi)]$ and $\tau = (J_a^*)^{-1}K_a^*(J_a^*)^{-1}$, provided the assumptions of \citet{jullum2017parametric} and \citet{daehlenaccurate} hold.
\end{theorem}

The crucial part of Theorem~\ref{thm:FIC theorem} is the fact that $\hat{b}^2 - n^{-1}(2bc + \kappa - \tau)$ is only $o(1/n)$ away from the true MSE of $g(\hat{\theta})$ in expected value. By Theorem~\ref{thm:HL limit marginal} the variance of $\hat{\theta}$ is of order $O(1/n)$. Hence, we need to approximate the MSE to a smaller order than $O(1/n)$. Otherwise the variance will be on the same scale as the error of the estimator of the MSE and might be greatly downplayed in the approximation.

In practice, the quantities $c$, $\kappa$ and $\tau$ are usually unknown. To approximate the MSE of the MHL estimator using Theorem~\ref{thm:FIC theorem}, we therefore need to first define estimators of these parameters and approximations to the expression $\hat{b}^2 - n^{-1}(2bc + \kappa - \tau)$. This is provided in the following definition.

\begin{definition}\label{def:HFIC}
	Let all quantities be as defined in Theorem~\ref{thm:FIC theorem}. We define the following approximation to the MSE of $g(\hat{\theta})$,
	\begin{align*}
		\text{HFIC} = \max\left\{\hat{b}^2 - n^{-1}(2\hat{b}\hat{c} + \hat{\kappa})\right\} + n^{-1}\hat{\tau},
	\end{align*}
	where $\hat{c}$ is either of the estimators given in \citet{daehlenaccurate}, $\hat{\tau} = \nabla g(\hat{\theta})^T(\hat{J_a^*})^{-1}\hat{K_a^*}(\hat{J_a^*})^{-1} \nabla g(\hat{\theta})$ and $\hat{\kappa}$ is given by $\hat{\tau} + n^{-1}\sum_{i=1}^n\phi(Y_i,\hat{\xi})^2 - 2\nabla g(\hat{\theta})^T(\hat{J_a^*})^{-1}\hat{V}$, with $\hat{J_a^*} = n^{-1}Hh_n(\hat{\theta})$, $\hat{K_a^*} =  n^{-1}\sum_{i=1}^n\nabla_{\hat\theta}\hat{\psi}(Y_i, \theta)\nabla_{\hat\theta}\hat{\psi}_{\theta}(Y_i, \hat{\theta})^T$ and $\hat{V} =  n^{-1}\sum_{i=1}^n\nabla_{\hat\theta}\hat{\psi}(Y_i, \theta)\phi(Y_i,\hat{\xi})^T$. Here $\hat{\psi}(Y_i, \theta)$ is equal to $\hat{\psi}(Y_i, \theta)$ with $\lambda(\theta)$ replaced by $\hat\lambda$, the solution to $0 = V_n(\lambda,\theta) =  \sum_{i=1}^nM_i(\theta)/[1+\lambda^TM_i(\theta)]$ and $\lambda'(\theta)$ by $-(S^1_n)^{-1}S_n^2$ where $S_n^1 = \partial_{\hat\lambda}V_n(\lambda,\hat\theta)$ and $S_n^2 =\partial_{\hat\theta }V_n(\theta,\hat\lambda)$.
\end{definition}
\begin{remark}
	In most applications the control and focus parameter will be the same, but there is no technical reason for this, and, if useful, $\mu$ and $\psi$ can be chosen to be different parameters. See Section~\ref{section:Wars}.
\end{remark}

The quantities $\hat{J_a^*}$, $\hat{K_a^*}$, $\hat{V}$ and $n^{-1}\sum_{i=1}^n\phi(Y_i,\hat{\xi})^2$ are consistent for $J_a^*$, $K_a^*$, $\textup{Var}[\phi(Y,\xi)]$ and $\textup{E}[\psi_\theta(Y,\theta_0)\phi(Y)]$ respectively under weak conditions. See e.g.~Lemma 16 of \citet{daehlenaccurate} for a set of sufficient conditions. When this is the case, $\hat{\tau}$ and $\hat{\kappa}$ converge in probability towards $\tau$ and $\kappa$ by the continuous mapping theorem and Theorem~\ref{thm:HL limit marginal}, making HFIC approximately unbiased for the mean squared error of the MHL estimator. Hence, to choose the balance parameter in a way minimizing the asymptotic mean squared error, we can simply compute the MHL estimator and corresponding HFIC score for some values of $a$ in a grid of points between zero and one. Afterwards, the balance parameter resulting in the smallest estimated mean squared error of the focus parameter can be chosen.

The above HFIC is only valid when plug-in estimators are not present. Thanks to Theorem \ref{thm:HL limit plug-in}, however, the criterion can also be defined in this context. We give the formulas below.

\begin{definition}
	Let all quantities be defined as in Theorem \ref{thm:HL limit plug-in} and assume the conditions of this theorem hold. Furthermore, let $\xi\in\mathbb{R}$ be some focus parameter taking the value $g(\theta)$ in the parametric model and $\hat\xi$ a consistent estimator of $\xi$ satisfying \eqref{eq:Influence function FIC} for some $\phi$. Define $\psi_{\textup{FIC}}(y,\theta,\lambda,P) =(m(y,\mu,P)^T/[1+\lambda^Tm(y,\mu,P)],\nabla_\theta \psi_P(y,\theta,P)^T, h(y,P)^T)^T$, and let $J\psi_{\textup{FIC}}(y,\theta,\lambda,P)$ denote its Jacobian matrix with respect to all parameters but $y$. Let $\hat A = n^{-1}\sum_{i=1}^nJ\psi_{\textup{FIC}}(Y_i,\hat\theta,\hat\lambda,\hat P)$ and $\hat B = n^{-1}\sum_{i=1}^n\psi_{\textup{FIC}}(Y_i,\hat\theta,\hat\lambda,\hat P)\psi_{\textup{FIC}}(Y_i,\hat\theta,\hat\lambda,\hat P)^T$ and $\hat C = \hat A^{-1}\hat B(\hat A^{-1})^T$. With this we define the HFIC as in Definition \ref{def:HFIC} with $\hat b = g[\hat\theta(\hat P)] - \hat\xi$, $\hat c$ equal to either of the estimators given in \citet{daehlenaccurate} and $\hat\tau = (\nabla g[\hat\theta(\hat P)],0)^T\hat C(\nabla g[\hat\theta(\hat P)]^T,0)^T$ and
	\begin{equation*}
		\hat\kappa = \hat\tau_P + n^{-1}\sum_{i=1}^n\phi(Y_i,\hat\xi)^2 - 2(\nabla_{\hat\theta(\hat P)}g(\theta)^T,0)\hat An^{-1}\sum_{i=1}^n\psi_{\textup{FIC}}[Y_i,\hat\theta(\hat P), \hat\lambda, \hat P]\phi(Y_i,\hat\xi)^T,
	\end{equation*} 
	where $\hat\lambda$ is the root of $\sum_{i=1}^nm[Y_i,\hat P,\hat\theta(\hat P)]/\{1 + \lambda^Tm[Y_i,\hat P,\hat\theta(\hat P)]\}$ and the matrix $(\nabla g[\hat\theta(\hat P)],0)^T$ is a block matrix of dimension $1\times (p+q+r)$ with zeros in the latter $q\times r$ components.
\end{definition}

By Lemma 5 in the supplementary material and the results of \citet{jullum2017parametric} and \citet{daehlenaccurate} the focused information criterion is only $O(1/n)$ away from the true MSE of $g[\hat\theta(\hat P)]$ in expected value. 

\section{Regression}\label{chapter:Regression}

So far, we have only worked with i.i.d.~data points. All of the results derived are, however, immediately applicable to regression settings if we are willing to treat the covariates as random variables. We will now illustrate how this can be done in the context of linear regression. These ideas were briefly visited in \citet{hjort2018hybrid}.

Let $Y\in\mathbb{R}$ be the response variable and $X$ the regressor, and assume that $Z_i = (Y_i,Y_i)$ for $i = 1,\ldots,n$ are i.i.d.~from some distribution $F$. Suppose furthermore that we wish to fit a linear regression model to these data. This corresponds to maximizing the function $\ell_n(\beta_0,\beta) = -\sum_{i=1}^n(Y_i-\beta_0 - \beta^TX_i)^2$. We will use $\ell_n$ as the likelihood function in the hybrid likelihood framework.

Suppose we are unsure if a subset of the covariates is relevant for predicting  the outcome of $Y$. In such a context, it could be beneficial to fit the model in a way giving additional weight to the components we know are important. This can be achieved with the hybrid likelihood framework. Let $\mathcal{S}$ be the set of indices corresponding to the predictors we believe are relevant for prediction, and for a general vector $v$ let $v_{\mathcal{S}}$ denote the components of $v$ with indices in $\mathcal{S}$. Then, the solution to $\textup{E}[m(Y,X,\beta)] = 0$ where $m(y,x,\beta_0, \beta_{\mathcal{S}}) = (y-\beta_0-\beta_{\mathcal{S}}^Tx_{\mathcal{S}})(1,x_{\mathcal{S}}^T)^T$, is the minimizer of $\textup{E}[(Y-\beta_0 - \beta_SX_{\mathcal{S}})^2]$. Hence, $m$ can be used to construct the empirical likelihood function for $(\beta_0,\beta_{\mathcal{S}}^T)^T$, the coefficients in the best linear fit when regressing $Y$ onto $X_{\mathcal{S}}$. 

Direct computations show that for the smaller regression of $Y$ onto $X_\mathcal{S}$, $(\beta_0(\beta),\beta_{\mathcal{S}}(\beta)^T)^T = [\textup{E} (X_{\mathcal{S}+1}X_{\mathcal{S}+1}^T)]^{-1}\textup{E}( X_{+1,\mathcal{S}}X_{+1})\beta = (\Sigma_{\mathcal{S}+1,\mathcal{S}+1} + \mu_{+1}\mu_{+1}^T)^{-1}(\Sigma_{\mathcal{S}+1,:} + \mu_{\mathcal{S}+1}\mu_{+1}^T)\beta$, where for any vector $v$, $v_{+1}$ indicates the vector $(1,v^T)$, $\Sigma_{\mathcal{S}+1,\mathcal{S}+1}$ is the covariance matrix of $X_{\mathcal{S}+1}$ and $\Sigma_{\mathcal{S}+1,:}$ is the covariance between $X_{\mathcal{S}+1}$ and the full component vector $X_{+1}$. Hence, for fixed $a$ we can define the hybrid log-likelihood function as $h_n(\beta_0,\beta,\mu,\Sigma) = (1-a)\ell_n(\beta_0,\beta) + a\log\textup{EL}_n[\beta_{\mathcal{S}}(\beta,\mu,\Sigma)]$, where $\textup{EL}_n$ is the empirical likelihood function of $\beta_\mathcal{S}$ constructed with the estimating function $m$. Maximizing $h_n$ for $a>0$ fits a linear regression model in all covariates, while giving additional attention to the fit of $Y$ regressed onto  $X_\mathcal{S}$.

In practice, estimation of $\mu$ and $\Sigma$ is rarely of interest. Because of this, it can be beneficial to replace them by $\hat\mu$, the empirical mean, and  $\hat\Sigma$, the empirical covariance, respectively rather than optimizing with respect to these variables. This corresponds to using the hybrid likelihood framework with plug-in estimators for $\mu$ and $\Sigma$. Indeed, let $h(y,x,\mu,\Sigma) = (x^T-\mu^T,(x-\mu)(x-\mu)^T-\Sigma)^T$, then $(\hat\mu,\hat\Sigma)$ satisfies $n^{-1}\sum_{i=1}^nh(Z_i,\hat\mu,\hat\Sigma) = o_p(1/\sqrt{n})$. Hence, inference of the hybrid likelihood estimator of $\beta$ using plug-in estimators for $\mu$ and $\Sigma$ can be made using Theorem \ref{thm:HL limit plug-in}. As a plug-in result is lacking from \citet{hjort2018hybrid}, this greatly improves the applicability of the preliminary ideas for regression which are discussed in \citet{hjort2018hybrid}.

\begin{remark}
	It is worth noting that when $\mathcal{S}$ consists of all indices, both $\ell_n$ and $\log\textup{EL}_n$ are maximized at $\hat{\beta}_{\textup{OLS}}$, the ordinary least square estimator of $\beta$. Hence, $\hat{\beta}_{\textup{OLS}}$ maximizes $h_n$ for all values of $a$ and is therefore equal to the MHL estimator regardless of the choice of balance parameter. This holds true, whether the parametric model used is specified correctly or not. When $\mathcal{S}\neq \{1,\ldots,d\}$, however, $\log \textup{EL}_n$ and $\ell_n$ do not have the same maxima and hence the MHL estimator is not constant as a function of $a$. 
\end{remark}

\section{Simulations}\label{section:Simulations}
\begin{figure}
	\centering
	\includegraphics[width=0.95\textwidth]{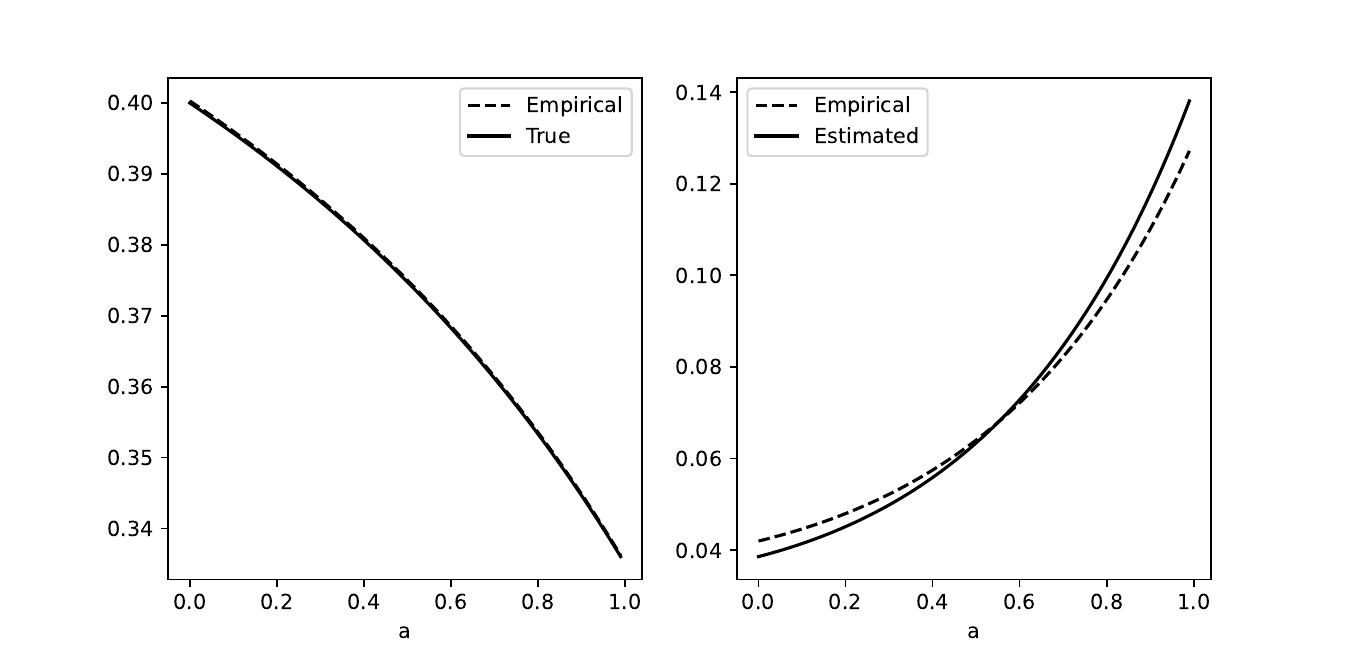}
	\caption{\label{fig:Simulation MHLE + var} The plot to the left shows the mean of $\hat p_1(a),\ldots,\hat p_B(a)$ for $a = 0,0.01,\ldots,0.99$ together with the theoretical limit. To the right, the empirical variance of $\sqrt{n}\hat p_1(a),\ldots,\sqrt{n}p_B(a)$ for each $a$ is displayed together with the mean of the estimator of the variance obtained using Definition \ref{def:HFIC}.}
\end{figure}

To make our results more trustworthy we performed a small simulation study, and did the following $B = 100$ times. We drew $n = 100$ data points from a Poisson distribution with rate parameter $1.5$ and fitted geometric distributions to the data with maximum hybrid likelihood for a selection of $a$-values. The geometric distribution has density $\log f_p(y) = y\log(1-p)+\log p$ for $y = 0,1,\ldots$. We used the probability of exceeding $1$ as the control parameter and $m(y,\mu) = I(Y\leq 2) - \mu$ as the estimating function in the empirical likelihood function. We will refer to the estimators as $\hat p_b(a)$ for $b = 1,\ldots,B$ and $a=0,0.01,\ldots,0.99$.  The mean of $\hat p(a)_1,\ldots,\hat p_B(a)$ for each of the selected $a$ values are shown in Figure \ref{fig:Simulation MHLE + var} together with the theoretical limit, the minimizer of $d_a(g,f_p)$ where $g$ is the pmf in the Pois$(2)$ distribution. From the figure, we notice that the mean is close to the theoretical limit. In each iteration we also computed the estimated variance of $\sqrt{n}\hat p(a)$ for each $a$ using the estimators defined in \ref{def:HFIC}. The mean of the $B$ estimator is shown in Figure \ref{fig:Simulation MHLE + var} together with the empirical variance of $\sqrt{n}\hat p_1(a),\ldots,\sqrt{n}p_B(a)$ for each $a$. Again, the theory and estimators seem to agree with the simulation results.

\begin{figure}
	\centering
	\includegraphics[width=0.95\textwidth]{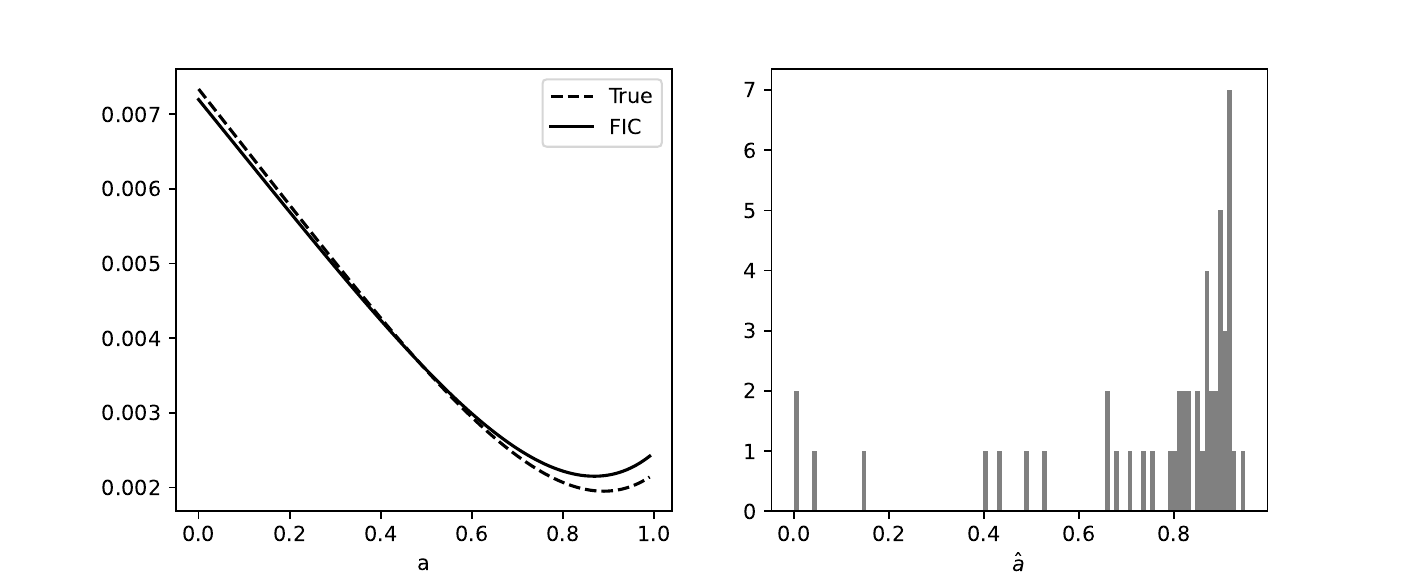}
	\caption{\label{fig:Simulation FIC + ahat} The plot to the left shows the mean of HFIC for $a = 0,0.01,\ldots,0.99$ together with the true MSE. To the right, a histogram of $\hat a_1,\ldots,\hat a_B$ is displayed.}
\end{figure}

We also computed the HFIC score for each $a$ in every iteration. We used the control parameter as the focus parameter and the empirical mean of $I(Y_i\leq 1)$ for $i = 1,\ldots, n$ as consistent estimator. This satisfies \eqref{eq:Influence function FIC} with $\phi(y,\xi) = I(y\leq 1) - \xi$. A plot of the mean of the $B$ HFIC curves is displayed in Figure \ref{fig:Simulation FIC + ahat} together with the observed mean squared error of the MHL estimator of the focus parameter. As before, the approximation seem to work very well. Since, the HFIC curves were computed in every iteration, we also registered the balance parameter minimizing the criterion in each case. A histogram of the $a$ chosen is shown in Figure \ref{fig:Simulation FIC + ahat}, and reveals that $a\neq0,0.99$ is often chosen, indicating that a true hybrid model is often preferred. Furthermore, the values of $a$ being most frequently chosen are those close to the minimizer of the MSE, indicating that HFIC is working as intended when it comes to choosing the balance parameter.

Lastly, we computed the MHL estimate of the focus parameter with the optimal $a$ chosen by HFIC in each iteration. This gave us estimates $\hat\mu_1(\hat a_1),\ldots,\hat\mu_B(\hat a_B)$. The MSE of these estimators were computed to be 0.003 which is smaller than the MSE of the ML estimators $\hat\mu_1(0),\ldots,\hat\mu_B(0)$ which was 0.007. This supports the claim that using hybrid likelihood in place of classic maximum likelihood can make results and estimators more precise.

\section{Application: Battle deaths in two hundred years of interstate wars}\label{section:Wars}
In his {\it The Better Angels of Our Nature}, \citet{pinker2012better} broadly argues that the world is changing for the better, and that the overall violence is on the decline. This claim quickly stirred up heated discussions in both the media and academic community, see e.g.~\citet{Guardian}, \citet{cirillo2016decline, cirillo2016statistical}. We will now illustrate how the hybrid likelihood framework can be applied to investigate certain aspects of these claims.

We extracted the number of battle deaths in the 95 most recent (and concluded) inter-state wars from the Correlates of War data set (\citet{CoW}). \citet*{cunen2020statistical} argue that the Korean War marks a change in the number of battle deaths, and hence, based on their analysis, we divided the data into two: older wars (the Korean war and all prior conflicts) and newer wars (the remaining data points). Our goal is to make inference about the difference between the number of battle deaths in newer and older wars. In \citet{daehlenaccurate} the authors investigate which models are suited to answer this question. They conclude that a shifted log-normal model outperforms the other options, including the Dagum model which is used in \citet*{cunen2020statistical}. Because of this, we will use shifted log-normal models for both data sets. 

We start by choosing a control parameter. Since $\mu\mapsto I(Y\leq \mu) - 1/2$ is discontinuous in $\mu$ and Theorem~\ref{thm:HL limit marginal} requires smooth estimating functions, we cannot control for the median directly. As a compromise we chose $m(y,\mu) = I(y\leq q_j)-\mu$ for $j=1,2$, where $q_1$ and $q_2$ are the empirical median number of battle deaths in older and newer wars, respectively. The theorems of this article can then be applied if $q_j$ are treated as fixed numbers and kept constant as $n\to\infty$.

\begin{figure}
	\centering
	\includegraphics[width=6in]{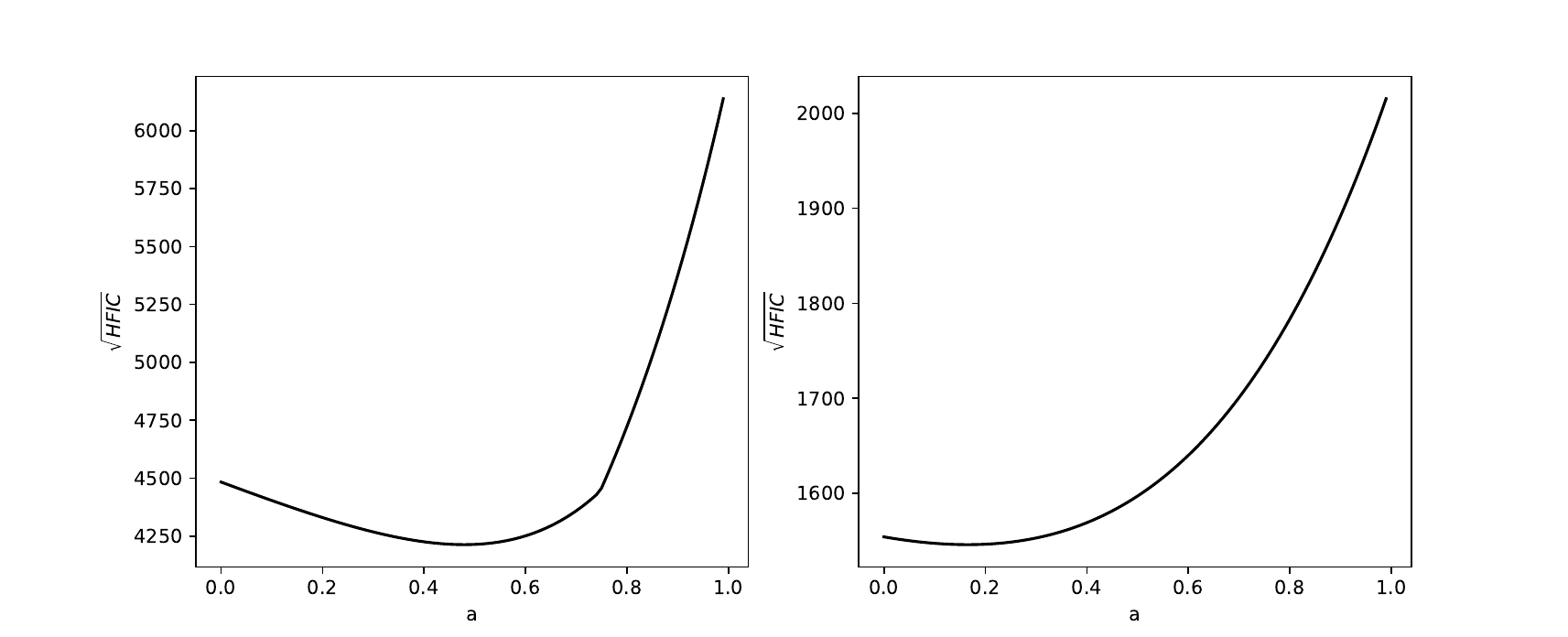}
	\caption{\label{fig:FICs Battledeaths}The figure displays root-HFIC scores as functions of $a$ (the balance parameter) for older and newer wars, in the left and right plots respectively. The minimum was attained at $a = 0.48$ for the older wars and $a = 0.17$ for the newer conflicts.}
\end{figure}

We chose separate balance parameters for each model.  In each case, this was done  by computing HFIC for each value of $a$ in a grid of points between $0$ and $1$. The median was used as the focus parameter for both models, and the empirical median in each data set as the consistent estimator $\hat{\xi}$ in Definition~\ref{def:HFIC}. Since the median satisfies \eqref{eq:Influence function FIC} with $\phi(x,\xi) = (I(Y\leq \xi) - 0.5)/f(\xi)$, where $f$ is the true density of the data, this is a valid choice. We estimated all quantities in HFIC as given in Definition~\ref{def:HFIC}, with $f(\mu)$ replaced by non-parametric kernel density estimates. For more details see \citet{daehlenaccurate}, in particular their Appendix A. The results are shown in Figure~\ref{fig:FICs Battledeaths}. We see that neither minima are attained at 0 or 1, i.e.~a genuine hybrid method is chosen for both data sets. The optimal balance values are $a=0.48$ for the older wars and $a=0.17$ for the newer conflicts. 

With all parameters set, we fitted the models to each data set using maximum hybrid likelihood. The resulting MHL estimators of the median for the two data sets were 8192 and 4750 for the older and newer wars respectively. In both cases this is closer to the empirical median in the data sets (11375 and 5240 respectively) than the maximum likelihood estimates (6799 and 4733). In addition, we estimated the standard deviation of the MHL estimator of the medians in the two populations. This was done by applying Theorem~\ref{thm:HL limit marginal} in conjunction with the delta method and the estimators defined in Definition~\ref{def:HFIC}. Estimated standard deviations were found to be 3488 and 1546 for the older and newer wars respectively. This is only marginally larger than the estimated standard deviations for the ML estimator for the older wars and more or less equal to the that of the the ML estimator for the newer wars (2746 and 1554 respectively). Hence, the added robustness of hybrid likelihood comes at a low price in this case.

We can now make inference about the difference in number of battle deaths between newer and older wars using hybrid likelihood. The MHL estimator of this parameter was found to be 3441, which, relatively speaking, is a lot larger than the corresponding ML estimator of 2066. We can compute variances for this quantity as well by applying the delta method and Theorem~\ref{thm:HL limit marginal}. We found the standard deviation for the MHL estimator to be 3815, while the ML estimator had a slightly smaller standard deviation of 3155. With these numbers we can construct confidence intervals. A 95\% confidence interval obtained with hybrid likelihood was computed to be $[-4036,10919]$. With maximum likelihood, we get a slightly smaller interval shifted to the left ($[-4135, 8232]$). In addition, we can construct hypothesis tests for testing if the median number of battle deaths is the same for older and newer wars, against the alternative that violence has gone down.  The p-value for testing this hypothesis using hybrid likelihood was 0.18, while the test based on maximum likelihood had a p-value of 0.26. These findings are in line with e.g.~\citet{clauset2018trends}.

The above analysis can be repeated for other quantiles than the median. To illustrate this, we performed similar computations as above for the difference between the upper quartiles in the two distribution. We leave out the details, but the p-value for testing if the difference is equal to zero was 0.03 for the hybrid likelihood framework. This is significant on a 5\%-level. The same is not true for the classical maximum likelihood method, which had a p-value of 0.06. Hence, the added robustness from the hybrid likelihood framework allows us to draw a more positive conclusion about the world. The difference is even clearer at higher quantile levels. Indeed, for each $p\in[0.50,0.85]$ we tested if the difference between the $p$-quantile for newer and older wars was zero against the alternative that violence has decreased. For almost all $p>0.65$ the null hypothesis is rejected. Hence, we can conclude that the larger wars have indeed become less deadly, which is consistent with the findings of \citet*{cunen2020statistical}.

The hybrid likelihood machinery may of course also be applied to other foci than the quantile difference, and with other control parameters. In particular, the two cumulative distribution functions, say $F_L(y)$ and $F_R(y)$, for intervals of $y$ of interest. The resulting estimates balance the parametric and nonparametric contributions. The difference is significantly non-zero for certain intervals of $y$.

It is worth noting that neither the ML estimator nor the MHL estimator are consistent for the true value of the parameter if the model is misspecified. The MHL estimator is, however, likely to aim for a quantity closer to the actual value of the parameter we are trying to estimate because of the added robustness introduced by the empirical likelihood part of $\textup{HL}_n$.

\section{Concluding remarks}
We have worked and significantly extended the hybrid likelihood framework of \citet{hjort2018hybrid}. Firstly, we derived alternative characterization of the empirical and hybrid likelihood functions. We also proved that the both the MEL and the MHL estimator are consistent for well-defined quantities and are asymptotically normal after proper scaling and centring. Our theory also allowed for the presence of plug-in parameters. Afterwards, we used the results to define estimators for the variance and MSE of the MHL estimator. We used this to define a focused information criterion and proposed to set the balance parameter to the minimizer of this criterion. Lastly, we showed how the hybrid likelihood framework could be applied to regression settings and illustrated the developed theory in an application. We will now briefly go through some limitations of our work and possible further research topics.

The perhaps most pressing issue with the theory developed in this article was discussed in Remark~\ref{remark:Bounded support} following Theorem \ref{thm:HL limit marginal}. If $m$ is the estimating function used in the construction of the empirical likelihood function, we require $m(y,\mu)$ to have bounded support for each fixed $\mu$ in a neighbourhood of $\mu_0$. Theoretically, this assumption is quite problematic, but as is shown in \citet[Thm.~1]{schennach2007point}, such a condition is necessary for the EL-framework to function outside of model conditions. In practice, however, it is less of a problem, as approximations using truncated distributions often give very good results. Take for instance, the situation where $m(y,\mu) = y-\mu$ and the data is normally distributed. Replacing each data point $Y$ by sign$(Y)\max\{K,|Y|\}$ for some large number $K$ makes all theorems go through while making virtually no difference to the analysis.

Another avenue for further research is to consider the case of non-smooth estimating functions $m$. The present theorems require $m$ to be differentiable. This excludes the case of quantiles as $m(y,\mu) = I(y\leq \mu)-p$ for $p\in[0,1]$ is not even continuous, let alone differentiable in $\mu$. The results do, however, hold true for smoothed approximations to $m$. Further investigation into the behaviour as the smoothed approximations get better and better might therefore lead to results covering the case of quantiles.

The regression setup discussed in Section~\ref{chapter:Regression}, required us to consider the covariates as random variables. This is more or less standard in literature concerning empirical likelihood (see e.g.~\citet{chen2009review}), but deviates somewhat from the classical approach in regression, where the covariates are often treated as fixed and non-stochastic. There is, however, ways of extending the empirical likelihood framework also to the non-i.i.d.~setting where the covariates are treated as non-random and fixed in regression settings, see e.g.~\citet{owen1991empirical}, \citet{wang1999empirical} or \citet{yan2015empirical}, and we therefore believe the framework of hybrid likelihood should be generalizable to this situation as well. In Lemma 3 in the supplementary material , we showed that the MHL estimator can be considered components of a Z-estimator. Combining this with classical results concerning solutions to generalized estimating equations (see e.g.~\citet[Chapter 5]{shao2003mathematical}) could lead to asymptotic theory and inference techniques for the MHLE estimator also in this non-i.i.d.~setup.

\section*{Supplementary Materials}

Proofs of Theorem \ref{thm:HL limit marginal} and \ref{thm:HL limit plug-in} can be found in the article's supplementary material.
\par
\section*{Acknowledgements}

This research has been partially funded by the Norwegian Research Council through the BigInsight Center for Research Driven Innovation (237718). The authors also gratefully acknowledge partial support from Centre for Advanced Study, at the Academy of Science and Letters, Oslo, in connection with the 2022-2023 project Stability and Change, led by H. Hegre and N.L. Hjort.
\par



  \bibliographystyle{chicago}      
  \bibliography{bib.bib}   

\vskip .65cm
\noindent
Ingrid Dæhlen
\vskip 2pt
\noindent
Department of mathematics, University of Oslo
\vskip 2pt
\noindent
E-mail: ingrdae@math.uio.no
\vskip 10pt
\noindent
Nils Lid Hjort
\vskip 2pt
\noindent
Department of mathematics, University of Oslo
\vskip 2pt
\noindent
E-mail: nils@math.uio.no

\newpage

\markboth{\hfill{\footnotesize\rm INGRID DÆHLEN AND NILS LID HJORT} \hfill}
{\hfill {\footnotesize\rm HYBRID LIKELIHOOD} \hfill}

\renewcommand{\thefootnote}{}
$\ $\par \fontsize{12}{14pt plus.8pt minus .6pt}\selectfont

\centerline{\large\bf  Supplementary material}
\noindent
The supplementary material contains proofs of Theorem 1 and 2 in the paper.
\par

\setcounter{theorem}{0}
\setcounter{section}{0}
\setcounter{equation}{0}
\def\theequation{S\arabic{section}.\arabic{equation}}
\def\thesection{S\arabic{section}}

\fontsize{12}{14pt plus.8pt minus .6pt}\selectfont

Here we will give formal proofs for the theorems in the article. The arguments utilize techniques for proofs concerning estimators fitted by two-stage procedures and are similar to those given by e.g. \citet*{imbens1995information} and \citet{schennach2007point}. Before we begin, however, we will state and prove a result concerning consistency of Z-estimators. The result is a modification of Theorem 5.42 in \citet{van2000asymptotic}.

\begin{lemma}\label{thm:Limit results Z-estimators}
	Let $Y_1,\ldots,Y_n\in\mathbb{R}^d$ be i.i.d. from some distribution $F$ and let $\psi\colon\mathbb{R}^d\times\Theta\to\mathbb{R}^q$ be a function with $\Theta\subseteq\mathbb{R}^p$. Assume there is a unique solution $\theta_n$ to the equation $0 = \Psi_n(\theta) = n^{-1}\sum_{i=1}^n\psi(Y_i,\theta)$ for all $n$ and a unique solution $\theta_0$ to the limit equation $0 = \Psi(\theta) = \textup{E}[\psi(Y,\theta)]$ when $Y\sim F$. Then $\theta_n$ converges in probability to $\theta_0$, provided
	\begin{itemize}
		\item [(C1)]
		the function $\Psi$ is continuously differentiable in a neighborhood of $\theta_0$ and $\Psi'(\theta_0)$ is invertible;
		\item [(C2)]
		there exists a neighborhood $\mathcal{N}$ of $\theta_0$ and an $F$-integrable function $p_0$ such that $\lVert\psi(y,\theta)\rVert\leq p_0(y)$ for all $\theta\in\mathcal{N}$ and $F$-almost all $y$.
	\end{itemize}
\end{lemma}
\begin{proof}
	Let $\Psi$ and $\Psi_n$ be defined as in the lemma. Inspection of the proof of Theorem 5.42 in \citet{van2000asymptotic} shows that the conditions of the theorem are slightly too strict and that the only thing needed for the arguments to go through is that $\Psi$ is continuously differentiable with non-singular derivative in a neighborhood of $\theta_0$ and that there exists a closed ball around $\theta_0$ inside which $\Psi_n$ converges uniformly to $\Psi$. The first condition is ensured by (C1) and continuity of the determinant function. The second condition follows from (C2) and the uniform law of large numbers. Hence, there exists a sequence of estimators converging in probability to $\theta_0$. By assumption, $\Psi_n(\theta) = 0$ has a unique solution $\theta_n$. Because of this, $\theta_n$ must be this consistent sequence, proving that $\theta_n$ converges in probability to $\theta_0$.
\end{proof}

\section{The maximum empirical likelihood estimator}\label{section:Proof EL}
We start by proving a limit result for the maximum empirical likelihood estimator (the MEL estimator). This will make the proofs concerning the maximum hybrid likelihood estimator more transparent.

\begin{lemma}\label{thm:Maximum empirical likelihood estimator}
	Let $Y_1,\ldots,Y_n\in\mathbb{R}^d$ be i.i.d. from some distribution $F$ and $m\colon\mathbb{R}^d\times\mathcal{M}\to\mathbb{R}^q$, where $\mathcal{M}\subseteq\mathbb{R}^s$ is an open set, be a function twice continuously differentiable in the last $s$ arguments. Let $\hat{\mu}$ be the maximizer of the empirical likelihood function and its unique critical point in $\mathcal{M}$. Furthermore, assume  $\hat{\lambda}$ is the unique solution to  
	\begin{align*}
		0 = n^{-1}\sum_{i=1}^n\frac{m(Y_i,\hat{\mu})}{1+\lambda^Tm(Y_i,\hat{\mu})},
	\end{align*}
	and $\hat{\alpha} = (\hat{\lambda}, \hat{\mu})$. Set   
	\begin{equation*}
		\psi_{\textup{EL}}(y,\alpha) = [1+\lambda^Tm(y,\mu)]^{-1}\left(m(y,\mu)^T,\,\lambda^T\partial m(y,\mu)\right)^T
	\end{equation*}
	with $\partial m$ as short hand for $\partial m/\partial\mu$. Lastly assume $\Gamma^{\textup{EL}}(\alpha) = \textup{E}\,\psi(Y,\alpha)$ has a unique root $\alpha_0 = (\lambda_0, \mu_0)$. Then $\hat\alpha$ is consistent for $\alpha_0$ and
	\begin{align}\label{eq:Limit of MELE}
		\hat{\alpha} = \alpha_0 - n^{-1}J\,\Gamma^{\textup{EL}}(\alpha_0)^{-1}\sum_{i=1}^n\psi_{\textup{EL}}(Y_i,\alpha_0)  + o_{\textup{pr}}(1/\sqrt{n})
	\end{align}
	if there exists a neighbourhood $\mathcal{N}$ on which
	\begin{itemize}
		\item [(1)]
		$\Gamma^{\textup{EL}}$ is continuously differentiable with non-singular Jacobian matrix at $\alpha_0$;
		\item [(2)]
		$\Pr(1+\lambda^Tm(Y,\mu) > L) = 1$ for some $L>0$ and all $(\lambda,\mu)\in\mathcal{N}$;
		\item [(3)]
		the norm of $m$ and all first and second order partial derivatives of $m$ are bounded by some functions $p_j$ for $j = 0,1,2$ respectively, with $\textup{E} [p_0(Y)^4],\,\textup{E} [p_1(Y)^4]<\infty$ and $\textup{E} [p_2(Y)^2]<\infty$.
	\end{itemize}
	In particular, $\sqrt{n}(\hat{\alpha}-\alpha_0)\overset{d}{\to}\text{N}[0,J\,\Gamma^{\textup{EL}}(\alpha_0)^{-1}K_{\textup{EL}}J\,\Gamma^{\textup{EL}}(\alpha_0)^{-1}]$,	with $K_{\textup{EL}} = \textup{Var}[\psi_{\textup{EL}}(Y,\alpha_0)]$, $Y\sim F$.
\end{lemma}

\begin{proof}
	Since $m$ is differentiable and $\mathcal{M}$ is open, the MEL estimator satisfies the equation
	\begin{align*}
		0 &= \frac{\partial}{\partial\mu}\log\textup{EL}_n(\mu) =  -\sum_{i=1}^n\frac{\lambda_n'(\mu)^Tm(Y_i,\mu) + \partial m(Y_i,\mu)^T\lambda_n(\mu)}{1+\lambda_n(\mu)^Tm(Y_i,\mu)} = -\sum_{i=1}^n\frac{\partial m(Y_i,\mu)^T\lambda_n(\mu)}{1+\lambda_n(\mu)^Tm(Y_i,\mu)},
	\end{align*}
	where the last equality follows from definition of $\lambda_n(\mu)$. Hence, $(\lambda_n(\hat{\mu}),\hat{\mu})$ is the solution to $0 = n^{-1}\sum_{i=1}^n\psi_{\textup{EL}}(Y_i,\alpha)$, where $\psi_{\textup{EL}}$ is as defined in the lemma. This shows that $\hat{\alpha}$ is a Z-estimator, and hence, standard results concerning such estimators can be used to study limiting properties of $\hat{\alpha}$. Furthermore, since the critical point of $\log\,\text{EL}_n$ was assumed to be unique. $\alpha_n$ is the unique root of $n^{-1}\sum_{i=1}^n\psi_{\textup{EL}}(Y_i,\alpha)$.
	
	We start by showing consistency towards $\alpha_0$. This is achieved by proving condition (C2) of Lemma~\ref{thm:Limit results Z-estimators}. (C1) holds by condition (1) of this lemma. By assumption (2) there exists $L>0$ such that $1+\lambda^TM(\mu)> L$ with probability 1 for all $\alpha\in\mathcal{N}$. Combining this with condition (3) ensures $\lVert\psi_{\textup{EL}}(y,\alpha)\rVert\leq C_1p_0(y) + C_2p_1(y)$ for some constants $C_1,C_2\in\mathbb{R}$. The right hand side of this inequality is a sum of $F$-integrable functions. This shows condition (C2) of Lemma~\ref{thm:Limit results Z-estimators}, and guarantees that $\hat{\alpha}$ converges in probability to $\alpha_0$.
	
	We will use Theorem 5.21 in \citet{van2000asymptotic} to show \eqref{eq:Limit of MELE}. We start with the Lipschitz condition. Direct computations show $\lVert\frac{\partial\psi}{\partial\alpha}(y,\alpha)\rVert\leq C_3p_0(y)^2 + C_4p_0(y)p_1(y) + C_5p_1(y) + C_6p_1(y)^2 + C_7p_2(y)$,	for some constants $C_j$ for $j=3,\ldots,7$. This is a consequence of condition (3) of the lemma. Since all functions on the right hand side of this equation are square integrable with respect to the measure $F$, the sum is square integrable as well. This in conjunction with the mean value theorem for multivariate functions, ensures the Lipschitz property of Theorem 5.21 in \citet{van2000asymptotic}. Furthermore,  $\textup{E}[\lVert\psi(Y,\alpha_0)\rVert^2]<\infty$ by a similar argument. The rest of the conditions of Theorem 5.21 in \citet{van2000asymptotic} hold by assumption. \eqref{eq:Limit of MELE} follows.
	
	Strictly speaking, we have only showed that the sequence of solutions to $n^{-1}\sum_{i=1}^n\psi(Y_i,\alpha)$ satisfies the properties of the lemma. Since, $\log\textup{EL}_n(\mu)$ is only equal to $-\sum_{i=1}^n\log[1+\lambda_n(\mu)^Tm(Y_i,\mu)]$ when $1+\lambda_n(\mu)^Tm(Y_i,\mu)>0$ for all $i=1,\ldots,n$, $\partial\log\textup{EL}_n(\hat{\mu})/\partial\mu = 0$ is not immediately equivalent to $n^{-1}\sum_{i=1}^n\psi(Y_i,\hat{\alpha}) = 0$. However, by construction the empirical likelihood function is always non-negative and equal to zero when $1+\lambda_n(\mu)^Tm(Y_i,\mu)\leq0$ for some $i=1,\ldots,n$. Hence, the maximizer of $-n^{-1}\sum_{i=1}^n\log[1+\lambda_n(\mu)^Tm(Y_i,\mu)]$ is equal to the maximizer of $\textup{EL}_n$ as long as it satisfies $1+\lambda_n(\mu)^Tm(Y_i,\mu)>0$ for all $i=1,\ldots,n$. Furthermore, since $1+\lambda^Tm(Y_i,\mu)>0$ for all $i=1,\ldots,n$ holds on $\mathcal{N}$ by assumption (2), $\hat{\alpha}$ will satisfy $1+\lambda_n(\mu)^Tm(Y_i,\mu)>0$ for all $i=1,\ldots,n$ and be the maximizer of $\log \textup{EL}_n(\mu)$ with probability tending to one.
\end{proof}

Next we show a result concerning the limiting properties of $\hat\mu$ alone.

\begin{theorem}\label{thm:EL limit marginal}
	Let all assumptions of Theorem \ref{thm:Maximum empirical likelihood estimator} hold true and all quantities be as defined there. Furthermore, assume that $\mu\to\textup{E}\log[1+\lambda(\mu)^Tm(Y,\mu)]$ where $\lambda(\mu)$ denotes the solution to $0 = \textup{E}\{m(Y,\mu)/[1+ \lambda^Tm(Y,\mu)]\}$ has a maximizer $\mu_0$. Then  the MEL estimator $\hat\mu$ is consistent for $\mu_0$, and 
	\begin{equation}\label{eq:Influence marginal MEL supplementary}
		\hat\mu = \mu_0 + n^{-1}J_{\textup{EL}}^{-1}\sum_{i=1}^n\psi_\mu(Y_i,\mu_0) + o_{\textup{pr}}(1/\sqrt{n}),
	\end{equation}
	with $J_{\textup{EL}} = H_{\mu_0}\textup{E}\{\log[1+\lambda^TM(\mu)]\}$ and $\psi_\mu(y,\mu) = -\nabla_\mu\log[1+\lambda(\mu)^Tm(y,\mu)]$, $\nabla_\mu$ being the gradient with respect to $\mu$.
\end{theorem}
\begin{proof}
	By Lemma~\ref{thm:Maximum empirical likelihood estimator}, $\hat{\mu}$ converges in probability to the last $s$ components of $\alpha_0$, the parameter solving $0 = \textup{E} [\psi_{\textup{EL}}(Y,\alpha)]$. By definition $\lambda(\mu)$ is the root of $ \textup{E} \{m(Y,\mu)/[1+\lambda^Tm(Y,\mu)]\}$, the first $q$ component functions in $\textup{E}[ \psi_{\textup{EL}}(Y,\alpha)]$. Hence, $\alpha_0 = (\lambda(\mu_0),\mu_0)$ where $\mu_0$ is the root of $\textup{E} \{\partial m(Y,\mu)^T\lambda(\mu)/[1+\lambda(\mu)^T m(Y,\mu)]\}$. Furthermore, the maximizer of $-\textup{E}\{ \log[1+\lambda(\mu)^Tm(Y, \mu)]\}$ satisfies $0 = \textup{E}\{[\lambda'(\mu)^Tm(Y,\mu) + \partial m(Y,\mu)^T\lambda(\mu)]/[1+\lambda(\mu)^Tm(Y,\mu)]\}$, since $\mathcal{M}$ is open and $m$ is differentiable. As $\lambda(\mu)$ is the root of $\textup{E}\{m(Y,\mu)/[1+$ $\lambda^Tm(Y,\mu)]\}$, the right hand side of this equation is equal to $\textup{E}\{\partial m(Y,\mu)^T\lambda(\mu)/[1+\lambda(\mu)^Tm(Y,\mu)]\}$. Since the root of $\Gamma^{\textup{EL}}(\alpha)$ is unique, $\mu_0$ must therefore be the the maximizer of $-\textup{E}\{ \log[1+\lambda(\mu)^Tm(Y,\mu)]\}$.
	
	For \eqref{eq:Influence marginal MEL supplementary}, let $(0,I)$ be the $(q+s)\times(q+s)$ block matrix where $I$ is the $s\times s$-identity matrix. Then, by \eqref{eq:Limit of MELE}, 
	\begin{align*}
		\hat{\mu} = (0,I)\hat{\alpha} = \mu_0 - (0,I)\Gamma^{\textup{EL}}(\alpha_0)^{-1}n^{-1}\sum_{i=1}^n\psi(Y_i,\alpha_0) + o_{\textup{pr}}(1/\sqrt{n}).
	\end{align*}
	Let $H_{xy}$ be the matrix obtained by first differentiating $\textup{E}\{\log[1+\lambda^TM(\mu)]\}$ with respect to $y$ and afterwards with respect to $x$ evaluated at $(\lambda_0,\mu_0)$. Then $J\,\Gamma^{\textup{EL}}(\alpha_0)$ can be written as the following block matrix 
	\begin{align*}
		J\,\Gamma^{\textup{EL}}(\alpha_0) = \begin{pmatrix}
			H_{\lambda\lambda}, &H_{\mu\lambda}\\
			H_{\lambda\mu}, &H_{\mu\mu}
		\end{pmatrix}.
	\end{align*}
	Furthermore, $\lambda'(\mu_0) = -H_{\lambda\lambda}^{-1}H_{\mu\lambda}$ by the implicit function theorem. Utilizing this and the formula for inverses of block matrices, we find $(0,I)\Gamma^{\textup{EL}}(\alpha_0) = (A^{-1}\lambda'(\mu_0)^T, A^{-1}) = A^{-1}(\lambda'(\mu_0)^T,I)$, where $A = H_{\mu\mu} - H_{\lambda\mu}H_{\lambda\lambda}^{-1}H_{\mu\lambda}$. Direct computation shows that $A = J_{\textup{EL}}$. Putting all of this together shows $-(0,I)\Gamma^{\textup{EL}}(\alpha_0)^{-1}\psi(y,\alpha_0) = J_{\textup{EL}}^{-1}(\lambda'(\mu_0)^T,I)\psi(y,\lambda_0,\mu_0) = J_{\textup{EL}}^{-1}\psi_{\mu}(y,\mu_0)$, proving \eqref{eq:Influence marginal MEL supplementary}.
\end{proof}

We also state and show a result for plug-in estimators.

\begin{lemma}\label{lemma:Result for plugin full}
	Let $Y_1,\ldots,Y_n$ be i.i.d. from some distribution $F$ and $m\colon\mathbb{R}^{d}\times\mathcal{M}\times\mathbb{R}^r\to\mathbb{R}^q$, where $\mathcal{M}\subseteq\mathbb{R}^s$ is open, be a function continuously differentiable in the last $s+r$ arguments. Assume furthermore, that $\hat P$ is some estimator satisfying $n^{-1}\sum_{i=1}^nh(Y_i,\hat P) = o_{\textup{pr}}(1\sqrt{n})$ and that the empirical likelihood function is constructed using $\hat P$ as a plug-in for $P$, i.e. using the estimating function $(y,\mu)\mapsto m(y,\mu,\hat P)$.  Furthermore, let $\hat\mu(\hat P)$ be the MEL estimator and unique critical point of $\log\textup{EL}_n$. Lastly, let $\hat\lambda$ be the unique solution to
	\begin{equation*}
		0 = n^{-1}\sum_{i=1}^n\frac{m(Y_i,\hat\mu,\hat P)}{1 + \lambda m(Y_i,\hat\mu,\hat P)}
	\end{equation*}
	and set
	\begin{equation*}
		\psi_{\textup{EL}}^P(y,\lambda,\mu,P) = [1 +\lambda^Tm(y,\mu,P)]^{-1}\left(m(y,\mu,P)^T,\lambda^T\partial_\mu m(y,\mu,P)\right)^T
	\end{equation*}
	with $\partial_\mu m$ being shorthand for $\partial m/\partial\mu$. Then, if $\alpha_{P,0} = (\mu_0,\lambda_0,P_0)$ is the unique solution to  $0 = \Gamma_P^{\textup{EL}}(\mu,\lambda,P) =  \left(\textup{E}[\psi_{EL}^P(Y,\alpha_P)]^T,\textup{E}[h(Y,P)]^T\right)^T$, the vector $\hat\alpha_P = (\hat\mu^T,\hat\lambda^T,\hat P^T)^T$ converges in probability to $\alpha_{P,0}$ and
	\begin{equation}
		\hat\alpha_P = \alpha_{P,0} - n^{-1}J\,\Gamma_P^{\textup{EL}}(\alpha_{P,0})^{-1}\sum_{i=1}^n\begin{pmatrix}
			\psi_{\textup{EL}}^P(Y_i,\alpha_{P,0})\\
			h(Y_i,P_0)
		\end{pmatrix} + o_p(1/\sqrt{n}),
	\end{equation}
	provided there exists a neighbourhood $\mathcal{N}$ of $\alpha_{0,P}$ on which:
	\begin{itemize}
		\item [(1)]
		the function $\Gamma_P^{\textup{EL}}$ is continuously differentiable with non-singular Jacobian matrix at $\alpha_{P,0}$;
		\item [(2)]
		$\Pr(1+\lambda^Tm(Y,\mu,P) > L) = 1$ for all $(\mu,\lambda,P)\in\mathcal{N}$ and some $L>0$;
		\item [(3)]
		the norm of $m$ and all first and second order partial derivatives of $m$ with respect to $\mu$ and $P$ are bounded by functions $p_j$ for $j = 0,1,2$ respectively, with $\textup{E} [p_0(Y)^4],\,\textup{E} [p_1(Y)^4]<\infty$ and $\textup{E} [p_2(y)^2]<\infty$.
		\item[(4)]
		the norm of $h(y,P)$ and all first order partial derivatives of $h$ with respect to $P$ is bounded by functions $q_1$ and $q_2$ not depending on $P$ and with $\textup{E}[q_1(Y)],\,\textup{E}[q_2(Y)^2]<\infty$.
	\end{itemize}
\end{lemma}
\begin{proof}
	By assumption, the plug-in estimator $\hat P$ satisfies $n^{-1}\sum_{i=1}^nh(Y_i,\hat P) = o_p(1/\sqrt{n})$. Hence, the vector $\hat\alpha_P$ is the solution to  
	\begin{equation*}
		0=n^{-1}\sum_{i=1}^n\begin{pmatrix}
			\psi_{\textup{EL}}^P(Y_i,\lambda,\mu,P)\\
			h(Y_i,P)
		\end{pmatrix} + \delta_n
	\end{equation*}
	where $\delta_n$ is a term of order $o_p(1/\sqrt{n})$ which does not depend on $\alpha_P$. The lemma now follows from arguments completely analogue to those given in the proof of Lemma \ref{thm:Maximum empirical likelihood estimator}.
\end{proof}

As with Lemma \ref{thm:Maximum empirical likelihood estimator}, the above result give asymptotic properties of the full vector $\hat\alpha_P$. In practice, however, $\hat\mu$ alone is of interest, and because of this, a marginal result for $\hat\mu$ alone is given below.

\begin{theorem}\label{thm:Maximum empirical likelihood estimator plug-in}
	Assume the conditions of Theorem \ref{thm:Maximum empirical likelihood estimator plug-in} hold true and let all quantities be as defined here. Furthermore, let $\lambda(\mu,P)$ and $P_0$ denote the solutions to $0 =\textup{E}\{m(Y,\mu,P)/[1+\lambda^Tm(Y,\mu,P)]\}$ and $\textup{E}h(Y,P) = 0$ respectively, and assume that the function  $\mu\mapsto\textup{E}\log[1+\lambda(\mu,P)^Tm(Y,\mu,P)]$ has a maximizer $\mu(P)$ for each $P$. Then the MEL estimator $\hat\mu(\hat P)$ is consistent for $\mu_0 = \mu(P_0)$ and
	\begin{equation}\label{eq:Influence plugin}
		\hat\mu(\hat P) = \hat\mu(P_0) + \mu'(P_0)(\hat P - P_0) + o_{\textup{pr}}(1\sqrt{n}),
	\end{equation}
	where $\hat\mu(P_0)$ is the MEL estimator when $P$ is fixed at $P_0$.
\end{theorem}

\begin{proof}
	Let $\alpha_{P,0} = (\lambda_0^T,\mu_0^T,P_0^T)^T$ be the root of E$[\Gamma_P^{\textup{EL}}(\alpha_P)]$ $=  (\textup{E}[\psi_{EL}^P(Y,\alpha_P)]^T,$ $\textup{E}[h(Y,P)]^T)^T$. Then, by definition, $P_0$ is the solution to E$[h(Y,P)]=0$ and $\lambda_0$ the solution to
	\begin{equation*}
		\textup{E}\left(\frac{m(Y,P_0,\mu_0)}{1+\lambda^Tm(Y,P_0,\mu_0)}\right) = 0.
	\end{equation*}
	Now let $\mu_1$ be the maximizer of $-\log[1+\lambda(\mu,P_0)^Tm(Y,P_0,\mu)]$ in $\mathcal{M}$. Then $\mu_1$ satisfies
	\begin{equation}
		0 = \textup{E}\left(\frac{\partial_{\mu_1}\lambda(\mu,P_0)^Tm(Y,P_0,\mu_1) + \lambda(\mu,P_0)^T\partial_{\mu_1}m(Y,P_0,\mu)}{1+\lambda(\mu,P_0)^Tm(Y,P_0,\mu_1)}\right).
	\end{equation}
	By definition of $\lambda(\mu_1,P)$ the second term equals zero, implying
	\begin{equation}
		0 = \textup{E}\left(\frac{ \lambda(\mu,P_0)^T\partial_{\mu_1}m(Y,P_0,\mu)}{1+\lambda(\mu,P_0)^Tm(Y,P_0,\mu_1)}\right),
	\end{equation}
	and ensuring that $\Gamma^{\textup{EL}}(\mu_1,\lambda(\mu_1,P_0),P_0) = 0$, but by assumption this equation has a unique root at $\alpha_{P,0}$. Hence, $\mu_1 = \mu_0$ and $\mu_0$ is the maximizer of $-\textup{E}\{\log[1+\lambda(\mu,P_0)^Tm(Y,P_0,\mu)]\}$.
	
	Secondly, consider the matrix $J\,\Gamma_P^{\textup{EL}}(\alpha)$. By definition of the Jacobian matrix, we can write $J\,\Gamma_P^{\textup{EL}}(\alpha)$ as the following block matrix:
	\begin{equation*}
		J\,\Gamma_P^{\textup{EL}}(\alpha) = \begin{pmatrix}
			J\,\Gamma^{\textup{EL}}(\lambda,\mu), &\partial_P\Gamma_P^{\textup{EL}}(\lambda,\mu,P)\\
			0, &\partial_P\textup{E}h(Y,P)
		\end{pmatrix},
	\end{equation*}
	where $J\,\Gamma^{\textup{EL}}(\lambda,\mu)$ is equal to $J\,\Gamma^{\textup{EL}}(\lambda,\mu)$ defined in Lemma \ref{lemma:Result for plugin full} where the estimating function used is$(y,\mu)\mapsto m(y,P,\mu)$ for each fixed $P$ and $\partial_P f$ is short-hand for $\partial f/\partial P$ for a general function $f$. By inversion of block matrices and condition (1) in Lemma \ref{thm:Maximum empirical likelihood estimator plug-in}, this matrix is continuously differentiable in a neighborhood of $(\lambda_0^T,\mu_0^T,P_0^T)^T$ with
	\begin{equation*}
		[J\,\Gamma_P^{\textup{EL}}(\alpha_{P,0})]^{-1} = \begin{pmatrix}
			[J\,\Gamma_P^{\textup{EL}}(\lambda,\mu)]^{-1}, &\alpha'(P_0)[\partial_P\textup{E}h(Y,P)]^{-1}\\
			0,&[\partial_P\textup{E}h(Y,P)]^{-1}
		\end{pmatrix}
	\end{equation*}
	non-singular. In the above $\alpha(P)$ is the root of $(\lambda,\mu)\mapsto\Gamma_P^{\textup{EL}}(\lambda,\mu,P)$ for each fixed $P$. Hence, by arguments similar to those given in the proof of Theorem 1,
	\begin{align*}
		&-\left(0\,I_s\,0\right)(J\,\Gamma_P^{\textup{EL}}(\alpha_{P,0}))^{-1}\psi_P(y,\alpha_{P,0}) =\\ &J_{\textup{EL}}^{-1}\psi_\mu^{P_0}(y,\mu_0) - \mu'(P_0)[\partial_P\textup{E}h(Y,P)]^{-1}h(y,P_0),
	\end{align*}
	where $\psi_\mu^{P_0}$ is equal to $\psi_\mu$ defined in Theorem 1 in the article when replacing $m$ by $(y,\mu)\mapsto m(y,\mu,P_0)$. Lastly, by condition (4) in Lemma \ref{thm:Maximum empirical likelihood estimator plug-in} and theorem 5.21 in \citet{van2000asymptotic}, $\hat P - P_0 = -[\partial_P\textup{E}h(Y,P)]^{-1}\sum_{i=1}^nh(Y_i,P_0) + o_p(1/\sqrt{n})$. Hence, $-\mu'(P_0)n^{-1}[\partial_P\textup{E}h(Y,P)]^{-1}$ $n^{-1}\sum_{i=1}^nh(Y_i,P_0) = \mu'(P_0)(\hat P - P_0) + o_p(1\sqrt{n})$. This concludes the proof.
\end{proof}

\section{Proving Theorem 1 in the article}
We first show a lemma.

\begin{lemma}\label{thm:Maximum hybrid likelihood estimator}
	Let $Y_1,\ldots,Y_n\in\mathbb{R}^d$ be i.i.d. from some distribution $F$ and $m\colon\mathbb{R}^{d+s}\to\mathbb{R}^q$ a function continuously differentiable in the last $s$ arguments. Let $f_{\theta}$ for $\theta\in\Theta\subseteq\mathbb{R}^p$, where $\Theta$ is an open set, be some parametric family with score function $u$ and $\mu\colon\mathbb{R}^p\to\mathbb{R}^s$ such that $\textup{E} \{m[Y, \mu(\theta)] \}= 0$ when $Y\sim f_{\theta}$. Fix $a\in [0,1)$ and assume that $\hat{\theta}$ is the maximizer and unique critical point of of $h_n$. Furthermore, let $\hat{\lambda}(\hat{\theta})$ be the unique solution to 	
	\begin{align*}
		0 = n^{-1}\sum_{i=1}^n\frac{m[Y_i,\mu(\hat{\theta})]}{1+\lambda^Tm[Y_i,\mu(\hat{\theta})]} = 0
	\end{align*} 
	Lastly, let $\beta = (\theta, \lambda)$, $\hat{\beta} = (\hat{\theta}, \hat{\lambda})$, and 
	\begin{align*}
		\psi_{\textup{HL}}(y,\beta) = 
		(1-a)\begin{pmatrix}
			0\\
			u(y,\theta)
		\end{pmatrix} - a\begin{pmatrix}
			I,&0\\
			0,&\mu'(\theta)^T
		\end{pmatrix}\psi_{\textup{EL}}[y,\lambda,\mu(\theta)].
	\end{align*}
	Assume that $\beta_0 = (\lambda_0, \mu_0)$ is the unique root of $\Gamma^{\textup{HL}}(\lambda, \mu) = \textup{E}[\psi(Y,\beta)]$, where $Y\sim F$. Then $\hat{\beta}$ is consistent for $\beta_0$ and 
	\begin{align}\label{eq:Limit of MHLE}
		\hat{\beta} = \beta_0 - n^{-1}J\,\Gamma^{\textup{HL}}(\beta_0)^{-1}\sum_{i=1}^n\psi_{\textup{HL}}(Y_i,\beta_0)  + o_{\textup{pr}}(1/\sqrt{n})
	\end{align}
	provided condition (2) and (3) of Lemma~\ref{thm:Maximum empirical likelihood estimator} hold true with $\alpha_0 = (\lambda_0^T,\mu(\theta_0)^T)^T$ and the following additional assumptions
	\begin{itemize}
		\item[(4)]
		the function $(\lambda,\theta)\mapsto\textup{E}\{(1-a)u(Y,\theta)^T -aM'(\theta)/[1+\lambda^TM(\theta)],M^T/[1+\lambda^TM(\theta)]\}$ is continuously differentiable in a neighbourhood of $\beta_0$ with non-singular derivative at $\beta_0$.
		\item[(5)]
		the function $\mu(\theta)$ is twice continuously differentiable in a neighbourhood of $\theta_0$.
		\item[(6)]
		in a neighbourhood of $\theta_0$, all components of $u$ and $\partial u/\partial\theta$ are bounded by some $F$-square integrable functions $q_1$ and $q_2$ respectively not depending on $\theta$.
	\end{itemize}
	In particular, $\sqrt{n}(\hat{\beta}-\beta_0)\overset{d}{\to}\text{N}(0,J\,\Gamma^{\textup{HL}}(\beta_0)^{-1}K_{\textup{HL}}J\,\Gamma^{\textup{HL}}(\beta_0)^{-1})$,	with $K_{\textup{HL}} = \textup{Var}[\psi_{\textup{HL}}(Y,\beta_0)]$, $Y\sim F$.
\end{lemma}
\begin{proof}
	Notice that
	\begin{align*}
		\frac{\partial h_n(\theta)}{\partial\theta} =
		& \sum_{i=1}^n\left\{(1-a)u(Y_i,\theta) - a\mu'(\theta)^T\frac{\partial m[Y_i,\mu(\theta)]^T\lambda_n(\theta) + \lambda_n'(\theta)^Tm[Y_i,\mu(\theta)]}{1+\lambda_n(\theta)^Tm[Y_i,\mu(\theta)]} \right\}\\
		=&\sum_{i=1}^n\left\{(1-a)u(Y_i,\theta) - a\mu'(\theta)^T\frac{\partial m[Y_i,\mu(\theta)]^T\lambda_n(\theta)}{1+\lambda_n(\theta)^Tm[Y_i,\mu(\theta)]}\right\},
	\end{align*}
	where the last equality follows from the definition of $\lambda_n(\theta)$ and $\partial m$ is short hand for $\partial m/\partial \mu$. Hence, since $\Theta$ is open and the maximizer of $h_n$ in $\Theta$ is the unique critical point of the function, $\hat{\beta}$ is the unique solution to $0 = n^{-1}\sum_{i=1}^n\psi_{\textup{HL}}(Y_i,\beta)$ where $\psi_{\textup{HL}}$ is as defined in the lemma. Because of this $\hat{\beta}$ is a Z-estimator and, just as for the MEL estimator, standard results can be used to study its limiting behavior.
	
	Let $m_2\colon\mathbb{R}^{d+p}\to\mathbb{R}^q$ be defined as $m_2(y,\theta) = m[y,\mu(\theta)]$. Since $\mu$ is twice continuously differentiable as a function of $\theta$, $m_2$ is twice continuously differentiable as a function of $\theta$. By the proof of Lemma~\ref{thm:Maximum empirical likelihood estimator} this ensures that there exists a neighborhood on which the norm of
	\begin{align}\label{eq:Helper in MHLE proof}
		\begin{pmatrix}
			I,&0\\
			0,&\mu'(\theta)^T
		\end{pmatrix}\psi_{\textup{EL}}[y,\lambda,\mu(\theta)] = \begin{pmatrix}
			-m_2(y,\theta)\\
			- \partial m_2(y,\theta)^T\lambda
		\end{pmatrix}\left[1+\lambda^Tm_2(y,\theta)\right]^{-1}
	\end{align}
	is bounded some square $F$-integrable function depending on neither $\lambda$ nor $\mu$. Combined with assumption (5), this proves that the norm of $\psi(y,\beta)$ is bounded by some square integrable function not depending on $\beta$. In particular, (C2) of Lemma~\ref{thm:Limit results Z-estimators} holds true. As (C1) holds by assumption (4), consistency of $\hat{\beta}$ towards $\beta_0$ follows.
	
	To show \eqref{eq:Limit of MHLE}, we will again use Theorem 5.21 in \citet{van2000asymptotic}. The only condition which remains to show is the Lipschitz condition. By the chain rule and proof of Lemma~\ref{thm:Maximum empirical likelihood estimator}, the norm of all partial derivatives of \eqref{eq:Helper in MHLE proof} are bounded by some square integrable function. Combining this with assumption (6) of this lemma and the mean value theorem, ensures that there exists a square integrable function $r$ such that $\lVert\psi_{\textup{HL}}(y,\beta_1)-\psi_{\textup{HL}}(y,\beta_2)\rVert\leq r(y)\lVert\beta_1-\beta_2\rVert$. This proves the Lipschitz condition of Theorem 5.21 in \citet{van2000asymptotic} and concludes the proof.
\end{proof}

\subsection{Proof of Theorem 1 in the article}
The proof of Theorem 1 in the article now follows directly from Lemma~\ref{thm:Maximum hybrid likelihood estimator} by performing computations similar to those in Section~\ref{section:Proof EL}. The case of $a=1$ follows from Theorem \ref{thm:EL limit marginal}.

\section{Proving Theorem 2 in the article}
We start with a lemma, yet again.

\begin{lemma}\label{thm:Maximum hybrid likelihood estimator plug-in}
	Let $Y_1,\ldots,Y_n\in\mathbb{R}^d$ be i.i.d. from some distribution $F$ and $m\colon\mathbb{R}^{d+s+r}\to\mathbb{R}^q$ a function continuously differentiable in the last $s+r$ arguments. Let $\hat P$ be some estimator satisfying $n^{-1}\sum_{i=1}^nh(Y_i,\hat P) = o_{p}(1/\sqrt{n})$ and assume that $(y,\mu)\mapsto m(y, \mu,\hat P)$ is used in construction of the empirical likelihood function $\textup{EL}_n$. Let $f_{\theta}$ for $\theta\in\Theta\subseteq\mathbb{R}^p$, where $\Theta$ is an open set, be some parametric family with score function $u$ and $\mu\colon\mathbb{R}^{p+r}\to\mathbb{R}^s$ such that $\textup{E} \{m[Y, \mu(\theta,P), P] \}= 0$ when $Y\sim f_{\theta}$. Fix $a\in [0,1)$ and assume that $\hat{\theta}(\hat P)$ is the maximizer of $h_n$ and unique critical value. Furthermore, let $\hat{\lambda}$ be the unique solution to 	
	\begin{align*}
		0 = n^{-1}\sum_{i=1}^n\frac{m[Y_i,\mu(\hat{\theta},\hat P),\hat P]}{1+\lambda^Tm[Y_i,\mu(\hat{\theta},\hat P),\hat P]} = 0
	\end{align*} 
	Lastly, let $\beta_P = (\theta, \lambda, P)$, $\hat{\beta} = (\hat{\theta}(\hat P), \hat{\lambda}, \hat P)$, and 
	\begin{align*}
		\psi_{\textup{HL}}^P(y,\beta) = 
		(1-a)\begin{pmatrix}
			0\\
			u(y,\theta)\\
			0
		\end{pmatrix} - a\begin{pmatrix}
			I,&0,&0\\
			0,&\mu'(\theta)^T,&0\\
			0,&0,&I
		\end{pmatrix}\begin{pmatrix}
			\psi_{\textup{EL}}^P[y,\lambda,\mu(\theta),P]\\
			h(y,P)
		\end{pmatrix}.
	\end{align*}
	Then if $\beta_{P,0} = (\theta_0,\lambda_0,P_0)$ is the unique root of $\Gamma^{\textup{HL}}_P(\lambda, \mu) = \textup{E}[\psi_\textup{HL}^P(Y,\beta)]$ where $Y\sim F$, $\hat{\beta}_P$ is consistent for $\beta_{P,0}$ and
	\begin{align}\label{eq:Limit of MHLE plugin}
		\hat{\beta} = \beta_0 - n^{-1}J\,\Gamma_P^{\textup{HL}}(\beta_0)^{-1}\sum_{i=1}^n\psi_{\textup{HL}}^P(Y_i,\beta_0)  + o_{\textup{pr}}(1/\sqrt{n})
	\end{align}
	provided condition (2-4) of Lemma~\ref{lemma:Result for plugin full} hold true with $\alpha_{P,0} = (\mu(\theta_0,P_0)^T, \lambda_0^T,P_0^T)^T$ and the following additional assumptions
	\begin{itemize}
		\item[(5)]
		The function $\Gamma^{\textup{HL}}_P$ is continuously differentiable in a neighbourhood of $\beta_{P,0}$ with non-singular derivative at $\beta_{P,0}$;
		\item[(6)]
		For all $P$ in an neighborhood of $P_0$, the function $\mu(\theta, P)$ is twice continuously differentiable with respect to $\theta$ and $P$ in a neighbourhood of $(\theta(P_0)^T, P_0^T)^T$.
		\item[(7)]
		There exists a neighborhood of $\theta_0$ on which all components of $u$ and $\partial u/\partial\theta$ are bounded by some $F$-square integrable functions $q_1$ and $q_2$ respectively not depending on $\theta$.
	\end{itemize}
\end{lemma}
\begin{proof}
	The result follows from arguments similar to those given in the proof of Lemma \ref{lemma:Result for plugin full} and Lemma \ref{thm:Maximum hybrid likelihood estimator} and is left out.
\end{proof}

\subsection{Proof of Theorem 2 in the article}
Theorem 2 now follows directly from the above lemma by arguing and performing computations similar to those in the proof of Theorem \ref{thm:Maximum empirical likelihood estimator plug-in}. The case when $a = 1$ follows from Theorem \ref{thm:Maximum empirical likelihood estimator plug-in}.

%

\end{document}